\documentclass{article}

\usepackage[english]{babel}
\usepackage{amsmath, amsxtra, latexsym,amscd, amsthm, amssymb,indentfirst}
\usepackage[mathscr]{eucal}
\usepackage{enumerate}
\usepackage{hyperref}
\usepackage{textcomp}
\usepackage{array, tabularx, longtable}%
\usepackage{multicol}
\usepackage{color}

\def\R{\mathbb R}
\def\N{\mathbb N}

\newcommand*{\rr}{{\rm{r}}}
\def\E{\mathbb E}

\def\int{{\rm int\,}}

\def\C{\mathcal C}
\def\B{\mathcal B}
\def\A{\mathcal A}
\def\F{\mathcal F}

\def\M{\mathcal M}

\newcommand*{\oa}{\overline{a}}

\def\det{{\rm det }}

\def\eps{\varepsilon}

\theoremstyle{plain}
\newtheorem{theorem}{Theorem}[section]
\theoremstyle{theorem}
\numberwithin{theorem}{section}
\newtheorem{corollary}[theorem]{Corollary}
\newtheorem{lemma}[theorem]{Lemma}
\newtheorem{proposition}[theorem]{Proposition}
\theoremstyle{definition}
\newtheorem{definition}[theorem]{Definition}

\newtheorem{example}[theorem]{Example}
\newtheorem{remark}[theorem]{Remark}

\newtheorem{notation}[theorem]{Notation}

\newtheorem{convention}[theorem]{Convention}
\begin{document}
	\title{Gaussian elimination for flexible systems of linear equations}
	\author{Nam Van Tran \\ Faculty of Applied Sciences, \\ HCMC University of Technology and Education, Vietnam
		\and Imme van den Berg\\ Research Center in Mathematics and Applications (CIMA), \\ 
		University of \'Evora, Portugal}
	\maketitle

	\begin{abstract}Flexible systems are linear systems of inclusions in which the elements of the coefficient matrix are external numbers in the sense of nonstandard analysis. External numbers represent real numbers with small, individual error terms. Using Gaussian elimination, a flexible system can be put into a row-echelon form with increasing error terms at the right-hand side. Then parameters are assigned to the error terms and the resulting system is solved by common methods of linear algebra. The solution set may have indeterminacy not only in terms of linear spaces, but also of modules. We determine maximal robustness for flexible systems.
	\end{abstract}
	
	\section{Introduction}

	We study systems of linear equations, with imprecisions in the coefficients and the right-hand side. We model the imprecisions asymptotically, however not functionally by $O(.)'s$ and $o(.)'s$ but instead by convex groups of nonstandard real numbers, called \emph{(scalar) neutrices}; 	
	we were inspired by Van der Corput's program for the Art of Neglecting \cite{Van der Corput}, with neutrices in the form of groups of functions, which are generalizations of the  $O(.)'s$ or $o(.)'s$ notations. A sum of a nonstandard real number and a neutrix is called an \emph{external number}. An external number can be seen as a real number with a small error term and captures the intrinsic vagueness of perturbations by the Sorites property of being invariant under some additions.
	
	A system of linear equations whose coefficients and right-hand side are given in terms of external numbers was called a \emph{flexible system} in \cite{Jus}. The Main Theorem of this article  presents a special form of Gaussian elimination applicable to all flexible systems, which is as effective as Gaussian elimination for real systems and leads to a solution in closed form, giving an explicit relation between the imprecisions of the system and the imprecisions of the solution. The method extends the Parameter method of \cite{Namopt} from non-singular systems to singular systems, and also generalizes the results of \cite{Jus} and \cite{NJI} on Cramer's Rule and Gauss-Jordan elimination for non-singular systems which are \emph{uniform}, i.e. all neutrices at the right-hand side are equal. 
	
	The central part of the solution method of a flexible system   applies Gaussian elimination in a careful form, always rearranging the system in order to be able to deal first with the smallest possible errors. In this way we obtain an equivalent system with real coefficient-matrix in \emph{increasing row echelon-form}, i.e. the coefficient matrix is in row echelon-form and the neutrices at the right-hand side corresponding to non-zero rows are increasing from top to down. The criterium for consistency of such a system is similar to the classical case: now the elements corresponding to zero rows in the coefficient matrix should be neutrices instead of zero's. If the system  is consistent, we apply the Parameter method and obtain an explicit solution in a form which is again similar to the solution of a classical non-determined system. Indeed, the solution set is the sum of a real vector and a  neutrix part, which is the solution of a homogeneous flexible system.   The neutrix part is the direct sum of a bounded neutrix and a linear subspace. The linear subspace is unique, but like the real support vector, the bounded neutrix is not unique. However it is the direct sum of scalar neutrices and has a well-defined dimension,  which is unique; we observe that neutrices are modules over $ \pounds $, the (external) set of limited real numbers. 
	
	The explicit formula improves the formula obtained in \cite{Namopt}, which still contemplated the intersection with the so-called feasibility space, a more-dimensional neutrix obtained from the neutrix parts of the coefficient matrix.
	
	We will apply the results on flexible systems to the problem of robustness for non-singular systems $ P|\B  $, where $ P $ is a real coefficient matrix and $ \B $ a vector with external numbers. This means we search for a matrix $ E\equiv (E_{ij}) $, where the $ E_{ij} $ are individual neutrices such that $ P|\B  $ and $ (P+E)|\B  $ are equivalent, i.e. have the same solution. 	We determine the maximal neutrices $ E_{ij} $ with this property, in the case that $ \det(P) $ is not too small.

	There are many approaches to study propagation of errors in linear systems, but to our knowledge in none of these settings straightforward Gaussian elimination has been applied to imprecise values in full generality. Among others methods to deal with imprecisions have been developed in the context of statistics and stochastics, fuzzy set theory, introduction of parameters, classical perturbation and error analysis, and interval calculus. 
	
	The latter methods are deterministic, hence conceptually are closer to our approach. We note that the calculation rules of external numbers (see Proposition \ref{defnam}) are the same as those for error analysis \cite{Taylor}. However, like in the case of interval calculus the formulation of general algebraic laws and advanced methods is restricted by complications due to the precise bounds of the error sets, for instance subdistributivity, intersection problems and loss of convexity; also upper bounds tend to be rapidly growing \cite{Alefeld}, \cite{Moore}, \cite{Gabrel}. These problems are still aggravated in the case of more variables \cite{Neumaier}.  As is the case of the present article, \cite{Assem} studies systems of the form \((A+\Delta A)(x+\Delta x)=b+\Delta b\). An upper bound of the relative error of the solution \(\dfrac{\Delta x}{x} \) with respect to \(\dfrac{\Delta b}{b} \) can be given with the help of the \emph{condition number} \({\rm{cond}}(A)=\|A\|\cdot \|A^{-1}\|\). However only examples of Gaussian elimination are given, with $ 2 $ and $ 3 $ variables.
	Disrespect of algebraic laws impeding a general theory of error propagation also occurs in the functional settings of parametric dependence \cite{Gal}, \cite{Gollmer} and classical assymptotics \cite{Bruijn}, \cite{Van der Corput}, in the latter case problems also arise from the difficulty to treat dependence of more variables and the lack of order.

	In the context of fuzzy set theory non-singular systems of any order were studied by among others B. Li and Y. Zhu in \cite{Liu}. In the case of a squared crisp matrix $ A $, a fuzzy right-hand side $ b $ and fuzzy variable $ x $ explicit solutions of the system $ Ax=b $ were given for two classes of distribution functions, of type exponential decay and piece-wise linear. The solution method involves matrix inversion for means and standard deviations. Fully fuzzy linear systems were studied by among others M. Dehghan, B. Hashemi and M. Ghatee \cite{Dehghan} using approximations by various well-known iterative methods. However to our knowledge there do not exist general results on error-propagation for Gaussian elimination in a fully fuzzy setting including singular systems. 
	
	This seems also true for sensitivity analysis based on statistics and stochastics, see e.g.,  \cite{Clifford}, \cite{Kovalenko}, \cite{Saltelli}. The study involves in one way or other the propagation of errors for operations on functions, and due to complexity, results in a general setting concern mainly properties of the solution, like mean, variance and bounds, than the solutions themselves.
	
	This article has the following structure. In Section \ref{secprel} we give some background on neutrices, external numbers and flexible systems. In Section \ref{secsol} we state the Main Theorem on the solution of flexible systems, describe the solution strategy and give an illustrative example. Section \ref{secalg} deals with the algebraic structure of neutrices in higher dimension. The various steps of the solution strategy are described in detail in Section \ref{secsolmet} and Section \ref{secprop} considers the algebraic structure of the solution. In Section \ref{secrank} we prove independence of rank when choosing real representative systems. The proof of the Main Theorem is completed in Section \ref{secproof}. In Section \ref{secrob} we present a model for robustness, and determine the maximal error allowed in each coefficient to not alter the solution of a flexible system. This section uses a result of Section \ref{secneg}, on the possibility to neglect rows of a system with small-enough errors.
	
	\section{Preliminaries}\label{secprel}
	
	\subsection{Scalar neutrices and external numbers}
	
	The article is written within the axiomatic form of nonstandard analysis $HST$ given by Kanovei and Reeken in \cite{Kanovei}. This is an extension of a bounded form of Internal Set Theory $IST$ of Nelson \cite{Nelson}, which in turn is an extension of common set theory $ ZFC $. To the language $ \{\in\} $ of $ ZFC $ a new predicate "standard" is added, denoted by "$ \rm{st}$". Formulas containing only the symbol  $ \{\in\} $ are called \emph{internal} and if they contain the symbol "$ \rm{st}$" they are called \emph{external}. Introductions to $ IST $ are contained in e.g. \cite{Dienerreeb}, \cite{Dienernsaip} or \cite{Lyantsekudryk}, and an introduction to a weak form of nonstandard analysis sufficient for a practical understanding of our approach is contained in \cite{DinisBerg}. An important tool is the principle of \emph{External induction} stating that induction is valid for all $ IST $-formulas over the standard natural numbers.

	The system $ IST $ distinguishes itself from Robinson's original model-theoretic approach \cite{Robinson}, by postulating that, next to the standard numbers, infinitesimals and infinitely large numbers already occur within the ordinary set of real numbers $\mathbb{R} $. A real number is \emph{limited} if it is bounded in absolute value by a standard natural number, and \emph{unlimited} if it is larger in absolute value than all limited numbers. Its reciprocals, together with $ 0 $, are called \emph{infinitesimal}. \emph{Appreciable} numbers are limited, but not infinitesimal.

	The notion "limited" refers to the predicate "standard", and the set $\pounds$ of all limited real numbers is an external subset of $ \R $ in the sense of $ HST $. Also the set of infinitesimals  $\oslash$, the set of positive unlimited numbers $ \not\hskip -0.14cm \infty $ and the set of positive appreciable numbers $ @ $ are external subsets of $ \R $. 
	
	The \emph{Minkowski operations} on subsets $ A,B $ of $ \R $ are defined pointwise. With respect to addition we have, with some abuse of language,
	\begin{equation*}\label{Minkowskisum}
		A+B=\{a+b|a\in A, b\in B\}.
	\end{equation*}
	The remaining algebraic operations on sets are defined similarly.
	
	\begin{definition}A \emph{(scalar) neutrix} is an additive convex subgroup of $\R$. An {\em external number}  is the Minkowski-sum of a real number and a neutrix. 
	\end{definition}
	Each external number  has the form $\alpha=a+A$, where $A$ is called the {\em neutrix part} of $\alpha$, denoted by $\rm{N}(\alpha)$, and $a\in \R$ is called a {\em representative}  of $\alpha$. We call $ \alpha $ {\em neutricial} if $ \alpha=\rm{N}(\alpha) $ and {\em zeroless} if $0\not\in \alpha$. The external class of all neutrices is denoted by $ \mathcal{N} $, this is not a proper external set in the sense of $ HST $, for "being a neutrix" amounts to an unbounded property \cite{DinisBergax}. The external class of all external numbers is denoted by $ \E $.
	
	The rules for addition, subtraction, multiplication and division of external numbers follow directly from the Minkowski operations. 
	
	\begin{proposition}\label{defnam}\rm Let $ a,b\in \R $, $A, B$ be neutrices and $\alpha=a+A, \beta=b+B$ be external numbers.
		\begin{enumerate}
			\item \label{cong} $\alpha\pm\beta=a\pm b+A+B$.		
			\item \label{nhan}  $\alpha \beta=ab+Ab+Ba+AB.$
			\item \label{1/alpha} If $\alpha$ is zeroless, $\dfrac{1}{\alpha}=\dfrac{1}{a}+\dfrac{A}{a^2}.$		
		\end{enumerate}
	\end{proposition}

	External neutrices are appropriate as a model for the Sorites property and orders of magnitude, for they are stable under some shifts, additions and multiplications. If an external number $\alpha=a+A$ is zeroless, one shows that its {\em relative imprecision} $ R(\alpha)\equiv \rm{N}(\alpha)/\alpha $ satisfies
	\begin{equation}\label{foutklein}
		{\rm{N}}(\alpha)=A/a\subseteq \oslash.
	\end{equation}
	So in combination with the intrinsic vagueness of the Sorites property, the neutrix $ A $ could be seen as a small error term for the real value $ a $ indeed. Observe that the algebraic rules of Proposition \ref{defnam} correspond to the rules of informal error analysis \cite{Taylor}. In particular we may recognize the property of neglecting the product of errors in the product rule given by Proposition \ref{defnam}.\ref{nhan}. Indeed, if $ \alpha $ or $ \beta$ is zeroless, by \eqref{foutklein} we have $ AB\subseteq Ab+Ba$, so we may neglect the neutrix product $ AB $. 
	
	A neutrix $ N$ is invariant under multiplication by appreciable numbers, i.e. $@N=N  $. An \emph{absorber} of $N $ is a real number $a $	such that $aN\subset N$ and an \emph{exploder} is a real number $b $	such that $bN\supset N$. So appreciable numbers are neither absorbers, nor exploders. Notions such as limited, infinitesimal, absorber and exploder may be extended in a natural way to external numbers.
	
	A neutrix $ N$ also satisfies $\pounds N=N  $, so from an algebraic point-of-view neutrices are modules over $ \pounds $. Division of neutrices is defined in terms of division of groups.
	
	\begin{definition}
		Let $A, B\in \mathcal{N}$. Then we define  
		$$ A:B=\{c\in \R\ |\ cB\subseteq A\}.$$
	\end{definition}

	An order relation for external numbers is given as follows.
	
	\begin{definition} Let $ \alpha,\beta \in \E$. We define 
		\begin{equation*}
			\alpha\leq\beta\Leftrightarrow\forall a\in \alpha \exists b\in \beta(a\leq b).
		\end{equation*}
		If $ \alpha \cap \beta =\emptyset $ and $ \alpha\leq\beta $, then $ \forall a\in \alpha \forall b\in \beta(a<b) $ and we write $ \alpha<\beta $.
	\end{definition}
	
	It is shown in \cite{Koudjeti Van den Berg} and \cite{DinisBerg} that the relation $ \leq  $ is an order relation indeed, which is compatible with the operations.
	If the neutrix $ A$ is contained in the neutrix $ B$, one has $ A\leq B $ and we say that $ B=\max(A,B) $. 
	
	The close relation to the real numbers and the group property of neutrices make practical calculations with external numbers quite straightforward. We always have subdistributivity and distributivity when multiplying with a real number, but in some cases related to subtraction distributivity does not hold. Necessary and sufficient conditions for distributivity are given in \cite[Theorem 5.6]{Dinis}.

	For more complete introductions to external numbers, including illustrative examples and lists of axioms, we refer to \cite{Koudjeti Van den Berg}, \cite{DinisBergax}, \cite{DinisBerg}.

	\subsection{Neutrices in higher dimension}
	Let $ n\geq 1 $ be standard. As in the one-dimensional case a set $ N\subseteq \R^n$ is called a \emph{neutrix} if it is a convex additive group. In analogy to external numbers we define external points as follows.
	\begin{definition}\label{defpoint}
		Let $ n\in \N $ be standard, $ p\in \R^{n} $ and $ N\subseteq\R^{n} $ be a neutrix. Then $ \xi=p+N\subseteq\R^{n} $ is called an \emph{external point}.
	\end{definition}
	
	Also in analogy to external numbers a \emph{representative point } $ p $ of an external point $  \xi $ is not uniquely determined, in contrast to the \emph{neutrix part} $ {\rm{N}}(\xi) =\{x'-x|x,x'\in \xi\}$.

	\begin{definition}\label{defboun}
		Let $ n\geq 1 $ be standard. A neutrix $ N\subseteq \R^{n} $ is called \emph{bounded} if there exists $ d\in \R,d>0 $ such that $ \lVert x\rVert <d$ for all $ x\in N $.
	\end{definition}
	
	\begin{definition}\label{deflinpar}
		Let $ n\geq 1 $ be standard and $ N\subseteq \R^{n} $ be a neutrix. The \emph{linear part} $ N_{(L)} $ of $ N $ is defined by
		\begin{equation}\label{NL}
			N_{(L)}=\cup\{S|S\subseteq N,S \mbox{ linear subspace of }\R^{n}\}.
		\end{equation} 
		A \emph{modular part} $ N_{(M)} $ of $ N $ is defined as a complement of $N_{(L)}$ in $ N $, i.e. it holds that $N_{(L)}\oplus N_{(M)}=N$. 
		Let $ x\in \R^{n} $ and $ \xi=x+N $ be an external point. Then we also call $ N_{(L)} $  the linear part of $ \xi $ and $ N_{(M)} $  a modular part of $ \xi $, and we write $ \xi_{(L)}=N_{(L)} $ and $ \xi_{(M)}=N_{(M)} $. 
	\end{definition}
	The linear part of a neutrix is uniquely defined, but a neutrix can have various modular parts. For instance $
	\R\times\oslash=\R\left( \begin{matrix}1\\0 
	\end{matrix}\right)+\oslash \left( \begin{matrix}0\\1 
	\end{matrix}\right)=\R\left( \begin{matrix}1\\0 
	\end{matrix}\right)+\oslash \left( \begin{matrix}1\\1 
	\end{matrix}\right)
	$.
	\begin{definition}\label{deforth}
		Let $ n\geq 1 $ be standard and $ N\subseteq \R^{n} $ be a neutrix. The \emph{ dimension} $ \dim(N) $ of $ N $ is defined by
		\begin{equation*}\label{fordimc}
			\dim(N)=\max\sharp\{ z|z\subseteq N, z \mbox{ linearly independent}\}.
		\end{equation*}
		
	\end{definition}
	A linear subspace $ V $ of $ \R^{n} $ is a particular case of a neutrix, and its dimension corresponds to the common dimension in the sense of linear algebra.
	
	\subsection{External vectors and matrices} \label{subsecvecmat}

	For a detailed account of vectors and matrices of external numbers we refer to \cite{Imme-nam-flexible matrix}. Here we recall some basic definitions and properties. 
	
	\begin{convention}
		From now on we always assume that $m,n\geq 1 $ are standard natural numbers.
	\end{convention}
	
	\begin{definition}\label{defextvec}
		Let $ A_1,\dots, A_n $ be neutrices. Then $A\equiv(A_1,\dots, A_n)^T$ is called a {\em neutrix vector}.
		
		Let $\B=(\beta_1,\dots, \beta_n)^T\in \E^n$, where $ \beta_i=b_{i}+B_{i} $ for $1\leq i\leq n$. Then $ \B $ is called an {\em external vector}. The vector $b=(b_1,\dots, b_n)^T$ is said to be a \emph{representative vector} of $\B$ and the neutrix vector $B=(B_1,\dots, B_n)^T$ is said to be the \emph{associated neutrix vector} of $\B$. 
		
		An external vector $\B=(\beta_1,\dots, \beta_n)^T\in \E^n$ can be identified with an external point with neutrix part in the form of a direct sum $ B_{1}e_{1}\oplus\dots \oplus B_{n}e_{n} $. However the notion of an external point is more general, for example, the external point $ \pounds \left( \begin{matrix}1\\1 
		\end{matrix}\right)\oplus\oslash\left( \begin{matrix}1\\-1 
		\end{matrix}\right)$ is not an external vector.
	\end{definition}
	\begin{definition}\label{defextmat}
		Let
		\begin{equation}\label{forextmat}
			\A=\begin{pmatrix}
				\alpha_{11} & \alpha_{12}& \cdots & \alpha_{1n}\\
				\vdots & \vdots&  \ddots & \vdots \\
				\alpha_{m1} & \alpha_{m2} & \cdots & \alpha_{mn}
			\end{pmatrix},
		\end{equation}
		where $\alpha_{ij}=a_{ij}+A_{ij}\in \E$ for $ 1\leq i\leq m, 1\leq j\leq n$. Then  $ \A $ is called an \emph{external matrix} and we use the common notation $\A=(\alpha_{ij})_{m\times n}$. For $ r< m $ we write $ \A^{rn} =(\alpha_{ij})_{1\leq i\leq r, 1\leq j\leq n}$. We denote by  $\M_{m, n}(\E)$ the class of all $m\times n$ external matrices. A matrix $\mathcal{A}\in \M_{m, n}(\E)$ is said to be \emph{neutricial} if all of its entries are neutrices, a special case is given by the zero matrix. We denote by  $\M_{m, n}(\R)$ the set of all $m\times n$ real matrices. With respect to  \eqref{forextmat} the matrix $P=(a_{ij})_{m\times n}\in \M_{m, n}(\mathbb{R})$ is called a \emph{representative matrix} and the matrix $A=(A_{ij})_{m\times n}$ the \emph{associated neutricial matrix}. If $ m=n $ we may write $ \M_{n}(\E) $ instead of $ \M_{n, n}(\E) $ and $ \M_{n}(\R) $ instead of $ \M_{n, n}(\R) $.
	\end{definition} 
	\begin{definition}\label{deflired}
		Let  $\B=(\beta_1,\dots, \beta_m)^T\in \E^m$ and $\A=(\alpha_{ij})_{m\times n}\in \M_{m, n}(\E)$.  
		We define 	$\left|\overline{\beta}\right|=\displaystyle\max_{1\leq i\leq m}|\beta_i|, $ 
		$ \underline{B}=\displaystyle\min_{1\leq i\leq m}B_i$, $\overline{B}=\displaystyle\max_{1\leq i\leq m}B_i$,  $ \overline{A_{i}}=\displaystyle\max_{\substack{1\leq j\leq n}} A_{ij} $ for $ 1\leq i\leq m $,
		$\left|\overline{\alpha}\right|=\displaystyle\max_{\substack{1\leq i\leq m\\1\leq j\leq n}} |\alpha_{ij}|$ and $
		\overline{A}=\displaystyle\max_{\substack{1\leq i\leq m\\1\leq j\leq n}} A_{ij}$. 	
		The external matrix $\mathcal{A}$ is said to be \emph{limited} if $ |\overline{\alpha}|\subset \pounds $ and \emph{reduced} if $ \overline{\alpha}=\alpha_{11} $ and $ \alpha_{11}=1+A_{11} $, with $ A_{11}\subseteq\oslash $, while all other entries have representatives which in absolute value are at most $1$. 
	\end{definition}
	By the last part of Definition \ref{deflired} a reduced external matrix always has  a reduced representative matrix.

	\begin{definition}
		For $ \A=(\alpha_{ij})_{m\times n},\A'=(\alpha'_{ij})_{m\times n}\in \M_{m,n}(\E)$ we write $\A\subseteq \A'$ if $\alpha_{ij}\subseteq \alpha'_{ij}$ for all $ i,j $ such that $ 1\leq i\leq m,  1\leq j\leq n$. 
	\end{definition}
	\subsection{Flexible systems} \label{subsecflex}
	
	Flexible systems were introduced in \cite{Jus} and studied also in \cite{Namopt} and \cite{NJI}. We recall the basic notions for flexible systems and introduce some new useful notions.

	\begin{definition}\label{defflex}
		Let 
		$ \A=  (\alpha_{ij})_{m\times n}\in \M_{m,n}(\E)$, $ x=(x_{1},\dots,x_{n})^{T}\in \R^{n} $ and  $ \B=(\beta_{1},\dots,\beta_{m})^{T} \in \E^{m} $. Then the set of linear inclusions  
		\begin{equation}\label{hpttqthamso}
			\left\{\begin{matrix}
				\alpha_{11} x_1+&\alpha_{12}x_2+&\cdots&+\alpha_{1n}x_n&\subseteq \beta_1\\
				\vdots &\vdots& \ddots&\vdots&\vdots\\
				\alpha_{m1} x_1+&\alpha_{m2}x_2+&\cdots&+\alpha_{mn}x_n&\subseteq \beta_m
			\end{matrix}\right.
		\end{equation}
		is called a \emph{flexible system} and denoted by $ \A x\subseteq\B $ or $ \A|\B $. 
		The  \emph{solution} of \eqref{hpttqthamso} is the set  $ \xi $ of all vectors $ x\in \R^{n} $ such that $ \A x\subseteq \B $. The solution is \emph{exact} if 	$ \A \xi=\B $.
	\end{definition} 
	
	The solution $ \pounds $ of the simple inclusion $ \oslash|\pounds $ is not exact, but as we shall see the solution of $ \A|\B $ is exact if $ \A $
	is a real matrix. 
	
	For flexible systems we will use throughout the notations of Definitions \ref{defextvec} and \ref{defextmat}.

	\begin{definition}\label{defsysA}
		The system $ \A|\B $ is called \emph{reduced} if $\mathcal{A}$ is reduced, \emph{limited} if $\A$ is limited, \emph{homogeneous} if $ \B $ is a neutrix vector, \emph{upper homogeneous} if $ \overline{\beta} $ is a neutrix and \emph{uniform} if the neutrices of the right-hand side $B_{i}\equiv B $ are all the same. When $ m=n $, the system is called {\em non-singular} if $\mathcal{A}$ is non-singular. 
	\end{definition}
	\begin{definition}\label{defeq}
		Let $ \A  \in \M_{m,n}(\E)$ and  $ \B, \B'\in \E^{m} $. Let $ p\in \N, p\geq 1 $ be standard and $ \A' \in \M_{p,n}(\E)$. The systems $ \A|\B $ and $ \A'|\B' $ are \emph{equivalent} if  $ \A x \subseteq \B\Leftrightarrow \A' x \subseteq \B'$ for all $ x\in \R^{n} $. Let $ H \in \M_{n}(\R)$ be a permutation matrix.
		We say that $ \A|\B $ and $ \A'|\B' $ are $ H $-\emph{equivalent} if $ \A x \subseteq \B\Leftrightarrow \A' y \subseteq \B'$ whenever $ x\in \R^{n} $ and $ y=Hx $.
	\end{definition}
	We will transform a general flexible system into a system with real coefficient matrix in \emph{increasing row-echelon form}, i.e. in row-echelon form, while the neutrices at the right-hand side of the non-singular part are increasing from above to below. It will be shown that such a system is equivalent to the original system up to renumbering the variables. 
	We recall first the notion of feasibility space, whose components  correspond to constraints for each individual variable \cite{Namopt}. We incorporate these constraints into the system, giving rise to the notion of integrated systems. Then we define the increasing row-echelon form, and finally we introduce some notation for the non-singular part of a system in increasing row-echelon form. The Main Theorem of Section \ref{secsol} affirms that the transformation can be carried out for any flexible system and the procedure is illustrated by  Example \ref{exfeas1}.
	
	\begin{definition}	Let $ \A\in \M_{m,n} (\E)$, $ \B\in \E^{m} $ and $ \A|\B $ be a flexible system.	For each $ j $ with $ 1\leq j\leq n $ we write
		\begin{equation*}
			F_{j}=\min_{1\leq i\leq m} B_{i}:A_{ij}.
		\end{equation*}
		The \emph{feasibility space} is defined by $  F\equiv\oplus_{j=1}^n F_j  e_j $.
	\end{definition}
	Note that some components of the feasibility space may be a vector space. In particular a component corresponding to a variable appearing only with real coefficients is equal to $ \R $, and a component corresponding to a variable appearing with a non-zero neutrix in a row with zero neutrix at the right-hand side is reduced to $ \{0\} $.

	\begin{definition}\label{deffeas}
		Let $ \A\in \M_{m,n} (\E)$, $ \B\in \E^{m} $ and $ \A|\B $ be a flexible system. Let $ F\equiv F_{1}e_{1}\oplus\cdots\oplus F_{n}e_{n} $ be the feasibility space corresponding to $ \A|\B $. Let $ k $ with $ 0\leq k\leq n $ be maximal such that $F_{j_{1}},\dots, F_{j_{k}}\subset\R $, with $ 1\leq j_{1}<\cdots<j_{k}\leq n $. Then we call $\F^{(c)}\equiv(F_{j_{1}},\dots, F_{j_{k}})^{T} $ the \emph{constraint}, $ k $ the \emph{constraint dimension}, and the  $ k\times n $-matrix $K=(d_{hl})_{1\leq h\leq k,1\leq l\leq n}$ defined by
		\begin{equation*}\label{Kron}
			d_{hl}=\left\lbrace \begin{array}{llc}
				1 && l=j_{h}\\
				0 &&\rm{else}
			\end{array}\right. 
		\end{equation*}
		the \emph{constraint matrix}.
	\end{definition}
	
	The matrix $ K $ is a sort of shifted identity matrix, with modified Kronecker symbols $ d_{hl} $, and indicates which variables do not range over the whole of $ \R $. Observe that  that $ K $ is uniquely determined and that if $ \A\in \M_{m,n} (\R) $, both $ \F^{(c)} $ and $ K $ are empty. 
	\begin{definition}\label{defintsys}
		Let $ \A\in \M_{m,n} (\E)$, $ \B\in \E^{m} $ and $ \A|\B $ be a flexible system. Let $ P $ be a representative matrix of \(\A\), $\F^{(c)} $ be the constraint having constraint dimension $ k $, and the  $ k\times n $-matrix $K\in \M_{k,n} (\R)$ be the constraint matrix. Then the system with real coefficient matrix \begin{equation}\label{realext}
			\left( \begin{matrix}
				P\\
				K
			\end{matrix} \left|  \begin{matrix}
				\mathcal B\\
				\mathcal F^{(c)}
			\end{matrix}\right. \right) 
		\end{equation}
		is called the \emph{associated integrated system} by $ P $.
	\end{definition}

	\begin{notation}
		Let $ Q=(q_{ij})_{m\times n}\in \M_{m,n} (\R)$ be of rank $ r\geq 1 $ and $ \mathcal{C} \in \E^{m} $. Then we write
		
		\begin{equation*}\label{Qr}
			Q^{(r)}\equiv\begin{pmatrix}
				q_{11}& \cdots & q_{1r}\\
				\vdots & \ddots & \vdots\\
				q_{r1} & \cdots & q_{rr}
			\end{pmatrix} ,
		\end{equation*}	
		the $ j^{th} $ column of $ Q $ by $ q_{j}^{T}$ and the first $ r $ elements of the $ j^{th} $ column of $ Q $ by $ (q^{r}_{j})^{T}$. We write $ \C=(\gamma_{1}, \dots, \gamma_{m})^T=(c_{1}+C_{1}, \dots, c_{m}+C_{m})^T=c+C$,  $ \C^{r}=(\gamma_{1}, \dots, \gamma_{r})^{T}$, $ \C^{m-r} =(\gamma_{r+1},\dots,\gamma_{m})^{T}$, $ c^{r}=(c_{1}, \dots, c_{r})^{T}$, and $ C^{r}=(C_{1}, \dots, C_{r})^{T}$. 
		If $ r<n $, for $ 1\leq i\leq n-r $ we denote the $i^{th}$ canonical unit vector in $\R^{n-r}$ by $e^{n-r}_i$.
	\end{notation}
	
	We now define the increasing row-echelon form, where for convenience we assume that the pivots are all on the principal diagonal.
	
	\begin{definition}\label{defincrow}
		Let $ Q\in \M_{m,n} (\R)$ be of rank $ r\geq 1 $, reduced, in row-echelon form and such that  $ q_{ii} =1$ if and only if $1\leq  i\leq r $, and $ \mathcal{C} \in \E^{m} $. We say that $ Q|\mathcal{C} $ is in \emph{increasing row-echelon form} if $ C_{1}\subseteq \cdots \subseteq C_{r } $.
	\end{definition}
	
	Observe that in Definition \ref{defincrow} we only require that the neutrix parts of right-hand side corresponding to non-zero rows of the coefficient matrix are  non-decreasing from top to down.
	
	\begin{definition}\label{assincrow}
		Let $ \A\in \M_{m,n} (\E)$, $ x\in \R^{n} $, $ \B\in \E^{m} $, and $ \A x\subseteq\B $ be a flexible system. Let $ P\in \M_{m,n} (\R)$ be a representative matrix of $ \A $.
		Assume the associated integrated system by $ P $	is $ H $-equivalent with a system $ Qy\in\C $  in increasing row-echelon form
		obtained by Gaussian elimination, where $ q=m+k$, $ Q\in M_{q,n} (\R)$,  $ \C\in\E^{q} $ and $ H\in \M_{n}(\R) $ is a permutation matrix. Then we call $ Qy\in\C $  a \emph{system in  increasing row-echelon form associated to $ \A x\subseteq\B $ by $ P $}. 
	\end{definition}

	\section{Main Theorem, solution set}\label{secsol}
	
	\begin{remark}\label{remboun}Let $ \A|\B $ be a flexible system, and $ x $ be a real vector. In some cases we apply a change of variables, and then we may make the variables explicit by writing, say, $ \A x\subseteq\B $, or $ \A x\in\B $ in case $ \A $ is real. We may still write the abbreviated form $ \A|\B $, if the variables are clear from the context, or if the symbols for the variables are not essential for understanding.
		
		We always assume that the neutrices at the right-hand side of a system $ \A|\B $ are different from $ \R $.
	\end{remark}

	The Main Theorem contains a general method to solve flexible systems $ \A|\B $, conditions for consistency and a closed form for the solution set.
	It gives additional information on the neutrix part of the solution set and the rank of the coefficient matrix of an associated integrated system.

	\begin{theorem}[Main Theorem]\label{maintheorem}

		Consider the system $ \A x\subseteq\B $, where $ \A \in \M_{m,n}(\E)$, $ x\in \R^{n} $ and $ \B=b+B \in \E^{m}$. Let $ P $ be a representative matrix of $ \A $ and $			\left( \begin{matrix}
			P\\
			K
		\end{matrix} \left|  \begin{matrix}
			\mathcal B\\
			\mathcal F^{(c)}
		\end{matrix}\right. \right) 	$ be an associated integrated system, where $ K \in \M_{k,n}(\R)$ is the constraint matrix and $ \F^{(c)} \in \E^{k}$ is the constraint, with $ k $ the constraint dimension. Let $ r\equiv\rr	\left( \begin{matrix}
			P\\
			K
		\end{matrix}\right) \geq 1 $.
		\begin{enumerate}
			
			\item \label{part1} There exists a permutation matrix  \(H\in \M_{n}(\R)\) such that  the system \(\A x\subseteq \mathcal{B}\)				
			is \(H\)-equivalent to a system \(Q y\in \mathcal{C}\) 
			in increasing row-echelon form which is associated to $ \A x\subseteq\B $ by $ P $, where $ Q\in \M_{m+k,n}(\R) $ has rank $ r $ and $ \mathcal{C}\equiv c+C \equiv(\gamma_{1},\dots,\gamma_{m+k})^{T} \in \E^{m+k}$, with $ \gamma_{i}=c_{i}+C_{i} $ for $ 1\leq i\leq m +k$. 	
			\item \label{part2} The system $ Qy\in\C $ is consistent if and only if $\gamma_{i}  $ is neutricial for $ r+1\leq i\leq m+k $; from now on, if we consider solutions and their properties, we will tacitly assume that $ Qy\in\C $, or equivalently $ \A x\subseteq\B $, is consistent.

			\item \label{part3} The solution  \(\zeta\) of $ Qy\in\C $ is exact, and given by 
			\begin{alignat}{2} \label{parasingular1}
			\zeta=\begin{pmatrix} \left(Q^{(r)}\right)^{-1}c^{r}\\
					0\end{pmatrix} +\sum\limits_{i=1}^r \begin{pmatrix} C_i\left(Q^{(r)}\right)^{-1}e_i^{r}\\
					0\end{pmatrix} +\sum\limits_{k=r+1}^n \R\begin{pmatrix} -\left(Q^{(r)}\right)^{-1} \left(q^{r}_{k}\right)^T\\
					e^{n-r}_{k-r}\end{pmatrix} .
			\end{alignat}  
			The matrix $ \left(Q^{(r)}\right)^{-1} $ is upper triangular. Moreover, the solution of the original system \(\A x\subseteq \B\) is given by \(\xi=H^{-1}\zeta\).			
			\item \label{part4} The linear part $ \zeta_{(L)} $ of $ \zeta $ is given by $ \zeta_{(L)}=\sum\limits_{k=r+1}^n \R\begin{pmatrix} -\left(Q^{(r)}\right)^{-1} (q^{r}_{k})^T\\
				e^{n-r}_{k-r}\end{pmatrix} $ and has dimension $ n-r $.
			\item \label{part5} The neutrix $ \zeta_{(M)}\equiv\sum\limits_{i=1}^r \begin{pmatrix} C_i\left(Q^{(r)}\right)^{-1}e_i^{r}\\ 
				0\end{pmatrix}  $ is a modular part  of $ \zeta $, and its dimension is equal to the number of non-zero neutrices of $ C^{r}$.	
			\item \label{part6} The rank of the coefficient matrix of the associated integrated system, the rank of the coefficient matrix of the associated integrated system, the neutrix part and 
			the linear part of the solution of $ \A x\subseteq \B $, and the dimension of its modular part do not depend on the choice of a representative matrix of $ \A $. 
		\end{enumerate}
	\end{theorem}

	The Main Theorem suggests the following solution method for flexible systems. To begin with we choose a representative matrix $ P $ of $ \A $,  write the constraints originating from the neutrix parts in the matrix form $ K|\F^{(c)} $ and join it to the system $ P|\B $ to obtain the associated integrated system $			\left(  \begin{matrix}
		P\\
		K
	\end{matrix} \left|  \begin{matrix}
		\mathcal B\\
		\mathcal F^{(c)}
	\end{matrix}\right. \right) 	$ with rank $ r $, say; it is shown in Subsection \ref{subsecint} that this is always possible. Then the integrated system is transformed into an equivalent system  $ Q\arrowvert \mathcal C $ in increasing row-echelon form; in Subsection \ref{subsecrow} it is shown that this can be done  by using a Gaussian elimination procedure involving, if necessary,
	interchanging columns of the coefficient matrix.
	Then we verify the condition for consistency, which simply amounts to verifying whether the components of $ \C $ of index bigger than $ r $ are neutricial.  
	In case of consistency we apply the parameter method of Theorem 4.3 of \cite{Namopt}. This means that parameters are assigned to the neutrices of $ \C $, then the system is solved by the usual means of linear algebra, which could be by repeated substitution, since the non-singular part of $ Q $ is upper triangular. Finally the closed form \eqref{parasingular1} is obtained by replacing the parameters in the solution formula by their range.
	
	Some elements of the solution procedure are not completely determined, in particular the choice of the representative matrix, and the choice of rows and columns in the Gaussian elemination process. The choices will influence the solution formula, in particular the support vector, the modular part and the basis of the linear part, in case the system is undetermined. However, Part \ref{part6} of Theorem \ref{maintheorem} ensures the invariance of the rank of the associated integrated matrix,  hence also of the rank of the associated matrix in increasing row-echelon form, of the linear part of the solution, and of the dimension of its modular part.
	
	The proof of Theorem \ref{maintheorem} consists of several steps. In Section \ref{secalg} we verify that  properties of linear algebra still hold for  neutrices. Section \ref{secsolmet} deals with the solution strategy sketched above. In Section \ref{secprop} we investigate the shape of the solution, and  in Section \ref{secrank} we prove invariance of ranks, when choosing representative matrices and vectors. In Section \ref{secproof} we put the results together and complete the proof. 
	
	Here we give an example illustrating the procedure sketched above.

	\begin{example}\label{exfeas1}
		Let $ \eps\simeq 0, \eps>0 $. Consider the flexible system 
		\begin{equation}\label{exsys1}	
			\left\{\begin{array}{rlllll}
				(-1+\eps\oslash)x_1+x_2 +(-0.2+\eps^{2}\pounds)x_3&\subseteq  2+\eps\pounds\\
				(1+\eps^{2}\pounds)x_1-x_2+(0.1+\eps^{2}\oslash)x_3&\subseteq 1+\eps\oslash\\
				(1+\oslash)x_1-x_2+(0.15+\eps\oslash)x_3&\subseteq -0.5+\oslash 
			\end{array}\right. .
		\end{equation}
		Then $ F_{1}=\min(\eps\pounds:\eps\oslash,\eps\oslash:\eps^{2}\pounds,\oslash:\oslash) =\pounds$ , $ F_{2}=\eps\oslash:\{0\}=\R $ and $ F_{3}=\min(\eps\pounds:\eps^{2}\pounds,\eps\oslash:\eps^{2}\oslash,\oslash:\varepsilon\oslash) =\pounds/\varepsilon$. Hence the feasibility space is given by \(\F =\pounds e_{1}\oplus \R e_{2}\oplus (\pounds/\varepsilon) e_{3}\).

		The constraint of the system \eqref{exsys1} is \(\F^{(c)}=(\pounds,\pounds/\varepsilon)^T\) and the constraint matrix is \(K=\begin{pmatrix}
			1&0&0\\
			0&0&1
		\end{pmatrix}.\) 
		
		We obtain a representative matrix $ P $ of the coefficient matrix of the system \eqref{exsys1} by neglecting the neutrix parts of the entries. Let  $ \B $ be the right-hand side, written in matrix form. Then the integrated system associated to the system \eqref{exsys1} becomes
		\begin{equation}\label{exsys2}
			\left(  \begin{matrix}
				P\\
				K
			\end{matrix} \left|  \begin{matrix}
				\mathcal B\\
				\mathcal F^{(c)}
			\end{matrix}\right. \right) =	\left( \begin{matrix}
				-1&1&-0.2	\\
				1&-1&0.1	\\
				1&-1&0.15	\\
				1&0&0\\
				0&0&1
			\end{matrix} \left|  \begin{matrix}
				2+\varepsilon\pounds\\
				1+\varepsilon\oslash\\
				-0.5+\oslash\\
				\pounds\\
				\pounds/\varepsilon
			\end{matrix}\right. \right).
		\end{equation}	
		
		We put now the system \eqref{exsys2} into increasing row-echelon form. The procedure asks first that all non-zero rows of the coefficient matrix are situated above the zero rows, which is trivially verified. Secondly each non-zero row should be reduced in such a way that some coefficient should equal to $ 1 $, while being maximal in absolute value. Again this is already verified.
		
		Then we interchange the first two rows since the neutrix at the right-hand side of the second row is smaller than the neutrix at the right-hand side of the first row. We get
		
		\begin{equation*}\label{exsys3}
			\left( \begin{matrix}
				1&-1&0.1	\\
				-1&1&-0.2	\\	
				1&-1&0.15	\\
				1&0&0\\
				0&0&1
			\end{matrix} \left|  \begin{matrix}
				1+\varepsilon\oslash\\
				2+\varepsilon\pounds\\
				-0.5+\oslash\\
				\pounds\\
				\pounds/\varepsilon
			\end{matrix}\right. \right).
		\end{equation*}
		Gaussian elimination of the first column leads to 
		\begin{equation*}\label{exsys4}
			\left( \begin{matrix}
				1&-1&0.1	\\
				0&0&-0.1	\\	
				0&0&0.05	\\
				0&1&-0.1\\
				0&0&1
			\end{matrix} \left|  \begin{matrix}
				1+\varepsilon\oslash\\
				3+\eps\pounds\\		
				-1.5+\oslash \\
				\pounds\\
				\pounds/\varepsilon
			\end{matrix}\right. \right).
		\end{equation*}
		We switch the second and third column, and obtain
		\begin{equation}\label{exsys5}
			\left(  \begin{matrix}
				1&0.1&-1	\\
				0&-0.1&	0\\	
				0&0.05&0	\\
				0&-0.1&1\\
				0&1&0
			\end{matrix} \left|  \begin{matrix}
				1+\varepsilon\oslash\\
				3+\eps\pounds\\		
				-1.5+\oslash \\
				\pounds\\
				\pounds/\varepsilon
			\end{matrix}\right. \right).
		\end{equation}
		We apply Gaussian elimination to the second column of \eqref{exsys5}, and obtain in a straightforward way
		\begin{equation}\label{exsys6}
			\left( \begin{matrix}
				1&0.1&-1\\
				0&1& 0\\		
				0&0&0\\ 
				0&0&1\\
				0&0&0
			\end{matrix} \left|  \begin{matrix}
				1+\varepsilon\oslash\\
				-30+\eps\pounds\\		
				\oslash \\
				\pounds\\
				\pounds/\varepsilon
			\end{matrix}\right. \right).
		\end{equation}
		Finally we interchange the third and the fourth row in \eqref{exsys6}, and get
		\begin{equation}\label{exsys7}
			\left( \begin{matrix}
				1&0.1&-1\\
				0&1& 0\\		
				0&0&1\\ 
				0&0&0\\
				0&0&0
			\end{matrix} \left|  \begin{matrix}
				1+\varepsilon\oslash\\
				-30+\eps\pounds\\		
				\pounds \\
				\oslash	\\
				\pounds/\varepsilon
			\end{matrix}\right. \right)\equiv Q|\mathcal{C}.
		\end{equation}
		The system $ Qy\subseteq\C $ of \eqref{exsys7} is in increasing row-echelon form indeed, with the neutrices $ (\varepsilon\oslash,\eps\pounds,\pounds) $ at the right-hand side corresponding to the non-singular part $ Q^{(3)}=			\left(\begin{matrix}
			1&0.1&-1\\
			0&1& 0\\		
			0&0&1
		\end{matrix}\right)  $ of the coefficient matrix $ Q $ increasing from above to below. 
		
		To solve the system, we ignore the last two rows, and let a parameter $ t_{1} $ range over $ \varepsilon\oslash $, $ t_{2} $ range over $ \varepsilon\pounds $ and $ t_{3} $ range over $ \pounds $. Then we get an ordinary upper triangular system, given by \begin{equation*}\label{exsys8}
			\left(\begin{matrix}
				1&0.1&-1\\
				0&1& 0\\		
				0&0&1\\ 	
			\end{matrix} \left|  \begin{matrix}
				1+t_{1}\\
				-30+t_{2}\\		
				t_{3}\\		
			\end{matrix}\right. \right).
		\end{equation*}
		We find $ (y_{1},y_{2},y_{3}) =(4+t_{1}-0.1t_{2}+t_{3},-30 +t_{2}, t_{3})$. Finally, noting that $ (x_{1},x_{2},x_{3})=(y_{1},y_{3},y_{2})  $, and substituting the parameters by their range we obtain the solution in vector form 
		
		\begin{equation}\label{solex}\xi\equiv
			\begin{pmatrix}
				\xi_1\\\xi_2\\\xi_3
			\end{pmatrix} =\begin{pmatrix}
				4 \\0\\-30
			\end{pmatrix}+ \eps\oslash \begin{pmatrix}
				1\\0\\0
			\end{pmatrix} + \eps\pounds \begin{pmatrix}
				-0.1\\0\\1
			\end{pmatrix}+ \pounds\begin{pmatrix}
				1\\1\\0
			\end{pmatrix}.
		\end{equation} 
		
		Geometrically, we could interpret the solution given by \eqref{solex} as a sort of affine space in the direction $ (1,1,0)^{T}$, truncated to $ \pounds $ and with support vector $ (4,0,-30)^{T}$, having a thin thickness $ \eps\pounds $ in the direction $ (-0.1,0,1)^{T}$ and still thinner thickness $ \eps\oslash $ in the direction $ (1,0,0)^{T}$.
	\end{example}
	\section{Algebraic properties of neutrices}\label{secalg}
	In analogy with systems of linear equations, the solution set of a flexible system $ \A|\B$ is the sum of a particular solution and the solution of a homogeneous inclusion. The latter is a neutrix. Theorem \ref{decNLnew} and \ref{themod} give additional information on neutrices in higher dimension. If a neutrix is a direct sum of a linear space $ V $ and a bounded neutrix $ W $, the former is necessarily equal to the linear part of the neutrix, and the dimension of $ W $ is uniquely determined. The latter is also true for a modular part, for as we will see, any modular part is a bounded neutrix. We will use these properties to prove Part \ref{part6} of the Main Theorem.
	\begin{proposition}Let $ \A \in \M_{m,n}(\E)$ and $ B \in \mathcal{N}^{m}$. Then		
		\begin{equation}\label{defN}
			N\equiv \{x\in \R^{n}|\A x\subseteq B \} 
		\end{equation}  
		is a neutrix.
	\end{proposition}
	
	\begin{proof}Let $ x,x'\in \R^{n} $.
		By subdistributivity $ \A(x-x')\subseteq\A x-\A x' \subseteq B+B=B$. Let $ 0\leq \lambda\leq 1 $. Again by subdistributivity $ \A(\lambda x+(1-\lambda)x')\subseteq \lambda Ax+(1-\lambda)Ax'\subseteq \lambda B+(1-\lambda)B=B$. We conclude that $ N $ is a convex group.
	\end{proof}
	
	\begin{theorem}\label{theoremsol}
		Let $ \A \in \M_{m,n}(\E)$ and $ \B=b+B \in \E^{m}$. Let $ \xi $ be the solution of the flexible system $ \A|\B$, and $ N $ be given by \eqref{defN}. If $ \xi $ is non-empty, it holds that $\xi=x+N $ for any $ x\in \xi $. Moreover $ {\rm{N}} (\xi)=N$.
	\end{theorem}
	\begin{proof}
		Let also $ y\in \xi$. Then 
		$  \A (y-x)\subseteq\A y-\A x\subseteq\B-\B=B$. So  $y-x\in N$, hence  $ y\in x+N $ and we derive that $ \xi\subseteq x+N $.
		
		Conversely, let
		$  y\in x+N$. Then $  y-x\in N$, so  $  \A (y-x)\subseteq B$, hence $ \A y =\A(x+ (y-x))\subseteq \A x+ \A(y-x))\subseteq \B+B=\B$. Hence  $y\in \xi$, which implies that $ x+N \subseteq \xi$.
		
		Combinining we obtain that $\xi=x+N $. Then $ {\rm{N}} (\xi)=N$.
		
	\end{proof}
	
	We show now that the linear part of a neutrix is uniquely determined, and any modular part is bounded, with uniquely determined dimension. 
	\begin{notation}\label{defVL}	Let $ n\geq 1 $ be standard and $ N\subseteq \R^{n} $ be a neutrix. We write 
		\begin{align*}
			N_{(\lambda)}&\equiv \cup\{\R v|\R v\subseteq N\}\\		N_{(\mu)}&\equiv \cup\{\R v|\R v\cap N\supset \{0\}\}\cup\{0\}.
		\end{align*}
	\end{notation}

	\begin{proposition}\label{lemlinpart}
		Let $ n\geq 1 $ be standard and $ N\subseteq \R^{n} $ be a neutrix.  Then $ N_{(L)} $ is a linear subspace of $\R^{n}  $ and
		\begin{equation}\label{forlemlin}
			N_{(L)}=N_{(\lambda)}.
		\end{equation}
	\end{proposition}
	\begin{proof}
		Firstly, we prove  \eqref{forlemlin}. The inclusion $ N_{(\lambda)}\subseteq N_{(L)}$ is obvious. Let $ v\in N_{(L)} $ and $V\subseteq N$ a linear subspace of $ \R^{n} $ such that $ v\in V $. Then $ \R v\subseteq N  $. So $ v\in N_{(\lambda)} $. Hence $ N_{(L)}\subseteq  N_{(\lambda)} $. Combining we obtain \eqref{forlemlin}.
		
		Secondly, we prove that $ N_{(L)} $ is a linear subspace of $\R^{n}  $. Clearly $ 0\in N_{(L)} $. Let $ v,w\in N_{(L)} $. Then, by the definition of \( N_{(L)}\), we have  $\R v\subseteq N  $ and $\R w\subseteq N $. In particular, $ v,w\in N $. Because $ N $ is a group, it holds that   $ v+w\in N $. Let $ t\in \R $. Then $ t(v+w)=tv+tw\in N+N=N $. Hence $\R (w+v)\subseteq N  $, meaning that $ v+w\in N_{(L)} $. We conclude that $ N_{(L)} $ is a linear subspace of $\R^{n}  $.
	\end{proof}
	\begin{proposition}\label{lemmodpart}
		Let $ n\geq 1 $ be standard and $ N\subseteq \R^{n} $ be a neutrix.  Then $ N_{(\mu)} $ is a linear subspace of $\R^{n}  $ and
		\begin{equation}\label{fordimmod}
			\dim(N_{(\mu)})=	\dim(N).
		\end{equation}
	\end{proposition}
	\begin{proof}
		Firstly, we prove that $ N_{(\mu)} $ is a linear subspace of $\R^{n}  $. Clearly $ 0\in N_{(\mu)} $. Let $ v,w\in N_{(\mu)} $.  Then there exist $ s\in \R, s\neq 0 $ such that $ sv,sw\in N $. Then $ s(v+w)=sv+sw\in N $. If $ w+v=0 $, clearly $ w+v\in N_{(\mu)} $. If $ w+v\neq 0 $, also $ s(w+v)\neq 0 $, so $\R (w+v)\cap N\supset\{0\}  $, which implies that $ v+w\in N_{(\mu)} $. Obviously $\R v\subseteq N_{(\mu)} $. We conclude that $ N_{(\mu)} $ is a linear subspace of $\R^{n}  $. 
		
		Secondly, let $ d= \dim(N), m=\dim(N_{(\mu)})$. If $  N=\{0\}$, also $  N_{(\mu)}=\{0\}$, and $ d=m=0 $. Assume now that $  N\supset\{0\}$, then also $  N_{(\mu)}\supset\{0\}$. There exist a linearly independent set $z\subseteq N$ such that $ d=\sharp z$. Because $ z\subset N_{(\mu)} $, we have $ d\leq m $. Conversely, there exists a linearly independent set of vectors  $w\equiv\{w_{1},\dots,w_{m}\}\subseteq N_{(\mu)}$. There exist $ s_{1},\dots,s_{m} \in \R\setminus\{0\}$ such that $ \{s_{1}w_{1},\dots,s_{m}w_{m}\}\subseteq N$. This implies that $ m\leq d $.
		
		Combining, we conclude that $ d=m $.
	\end{proof}
	
	\begin{theorem}\label{decNLnew}
		Let $ n\geq 1 $ be standard and $ N\subseteq \R^{n} $ be a neutrix. Let 
		$ V $ be a linear subspace of $ \R^{n} $ and $ W \subseteq \R^{n}$  be a bounded neutrix and  such that $ N=V\oplus W $. Then 
		\begin{enumerate}
			\item \label{decNLnew1} $ V=N_{(L)} $.
			\item \label{decNLnew2} $ \dim(W)=\dim(N)-\dim (N_{(L)} )$
		\end{enumerate}
	\end{theorem}
	\begin{proof}

		\ref{decNLnew1}. If $  V=\{0\} $, or \(V=\R^n\) it is easy to see that \(N_{(L)}=V\). In the remaining case $ \{0\}\subset V \subset \R^{n}$.
		
		By \eqref{NL} it holds that $ V\subseteq N_{(L)} $. Conversely, we show first that \(N_{(\lambda)} \subseteq V\). Observe that since \(W\) is bounded, there exists \(b\in \R,b>0\) such that for all $x\in W  $
		\begin{equation}\label{wbb}
			\|x\|\leq b.
		\end{equation}
		
		Let $ v\in N_{(\lambda)}, v\neq 0 $.  Suppose that $ v\notin V $. Let $ k= \dim (V)$. Then  $ 1\leq k\leq n-1$. Let  $ \{v_{1},\dots\,v_{k}\} $ be an orthonormal basis of $ V $. Because $ V $ is a linear subspace of $ \R^{n} $, it holds that $ \{v_{1},\dots,v_{k},v\} $ is linearly independent. Let $U$ be the linear subspace of $ \R^{n} $ spanned by $ \{v_{1},\dots,v_{k},v\} $. Then $ U\subseteq N+N=N $. Applying the Gram-Schmidt orthogonalization procedure, we find a vector $ v_{k+1} $ such that $ \{v_{1},\dots,v_{k},v_{k+1}\} $ is an orthonormal set of vectors in $ U $. Then also $ \R v_{k+1} \subseteq U\subseteq N$.  		
		We may complete $ \{v_{1},\dots\,v_{k+1}\} $ to an orthonormal basis $ \{v_{1},\dots\,v_{n}\} $ of $ \R^{n} $.
		
		Let $ t\in \R, t>b$. Then  
		\begin{equation}\label{toutn}
			\langle tv_{k+1}, v_{k+1}\rangle >b.
		\end{equation} 
		It follows from the fact \(\R v_{k+1}\subseteq N\) that \(tv_{k+1}\in N=V\oplus W\). So \(tv_{k+1}=y+w\) where \(y=y_1v_1+\cdots + y_kv_k \in V\) and \(w=w_1v_1+\cdots+w_nv_n\in W\), with \(y_1,\dots,y_k,w_{1},\dots,w_{n}\in \R\). 
		Hence
		\begin{align*}
			(w_1+y_1)v_1+\cdots +(w_k+y_k)v_k+ w_{k+1}v_{k+1}+\cdots +w_n v_n&=tv_{k+1}\\
			&=\langle tv_{k+1},v_{k+1}\rangle v_{k+1}.
		\end{align*}
		Because of the uniqueness of the representation of the vector \(tv_{k+1}\) in the basis \(\{v_1, \dots, v_n\}\) one has \(w_{k+1}=\langle tv_{k+1},v_{k+1}\rangle\). Then $ |w_{k+1}|>b $ by \eqref{toutn}, while 	$ |w_{k+1}|\leq \|w\|\leq b $ by \eqref{wbb}, a contradiction. Hence $ v\in V $, and we derive that $ N_{(\lambda)}\subseteq V $. Because $ N_{(\lambda)}=N_{(L)} $, it holds that $ N_{(L)} \subseteq V  $. 
		
		Combining, we conclude that $ N_{(L)}=  V  $.
		
		\ref{decNLnew2}. Because $ V\cap W=\{0\} $, also $ V\cap W_{(\mu)}=\{0\} $. Hence $ V\oplus W_{(\mu)}=N_{(\mu)}$, and being all linear spaces, $ \dim(V)+\dim\left( W_{(\mu)}\right)= \dim(N_{(\mu)})$. Then it follows from \eqref{fordimmod} that
		\begin{equation*}
			\dim(W)=\dim\left(W_{(\mu)}\right)=\dim\left(N_{(\mu)}\right)-\dim(V)=\dim(N)-\dim \left(N_{(L)}\right).
		\end{equation*}
	\end{proof}
	We recall now a definition and a theorem of \cite{BergAPAL}.
	\begin{definition}\label{delength}
		Let $ n\geq 1 $ be standard and $ N\subseteq \R^{n} $ be a neutrix. Then the neutrix
		\begin{equation*}\label{forlength}
			\varLambda(N)\equiv\{\lambda\in\R|\exists u\in\R^{n} \| u\|=1,\lambda u\subseteq N\}
		\end{equation*}
		is called the \emph{length} of $ N $.
	\end{definition}
	
	\begin{theorem}{\rm{\cite[Theorem 5.2]{BergAPAL}}}\label{thelength}
		Let $ n\geq 1 $ be standard and $ N\subseteq \R^{n} $ be a neutrix with lenght $ \varLambda $. Then there exists a unit vector $ u \in \R^{n}$ such that $ \R u \cap N=\varLambda u $.
	\end{theorem}
	\begin{theorem}\label{themod}
		Let $ n\geq 1 $ be standard and $ N\subseteq \R^{n} $ be a neutrix. Let 
		$ M $ be modular part of $ N $. Then 
		\begin{enumerate}
			\item \label{themod1} $ M $ is a bounded neutrix.
			\item \label{themod2} $ \dim(M)=\dim(N)-\dim (N_{(L)} )$.
		\end{enumerate}
		\begin{proof}
			As for Part \ref{themod1}, suppose $ \varLambda(M)=\R $. By Theorem \ref{thelength} there exists $ u\in \R^{n} $ such that $  \R u \cap M=\R u  $, i.e. $ \R u\subseteq M $. Hence $ M\cap N_{(L)}\supset \{0\} $, a contradiction. Hence $ \varLambda(M)\subset\R $. This implies that  $ M \subseteq (\varLambda(M))^{n}$ is bounded. Then Part \ref{themod2} follows from Theorem \ref{decNLnew}.\ref{decNLnew2}.
		\end{proof}
	\end{theorem}
	
	\section{Solution strategy}\label{secsolmet}

	\subsection{Integrated system}\label{subsecint}

	Let $ \A|\B $ be a flexible system. Theorem \ref{thdec} states that an associated integrated system is equivalent with the original system. In a sense it is a reformulation of \cite[Th.3.3]{Namopt}.

	\begin{theorem}\label{thdec}
		Let  \(\A=(\alpha_{ij})_{m\times n}\in \M_{m,n}(\E)\) with \(\alpha_{ij}=a_{ij}+A_{ij}\) for  $1\leq i\leq m, 1\leq j\leq n$, and $ \B=(b_1,\dots, b_n)^T+(B_1,\dots, B_n)^T \in \E^{m}$. Let $ P=(a_{ij})_{m\times n}\in \M_{m,n}(\R) $ be a representative matrix of \(\A\). Let $ k\in \N $ be the constraint dimension, $ \F^{(c)}\in \mathcal{N}^{k} $ be the constraint and  $K\in \M_{k,n} (\R)$ be the constraint matrix. Then the associated integrated system $			\left( \begin{matrix}
			P\\
			K
		\end{matrix} \left|  \begin{matrix}
			\mathcal B\\
			\mathcal F^{(c)}
		\end{matrix}\right. \right)	$ is equivalent to $ \A|\B $. 		
	\end{theorem}
	\begin{proof}  
		A vector  $x=(x_1, \dots, x_n)^T\in \R^n$ is a solution of the system $ \A|\B $ if and only if 
		$$\left\{ \begin{array}{ccccccccclllllll}
			(a_{11}+A_{11}) x_1& +&\cdots&+&(a_{1n}+A_{1n})x_n&\subseteq &  b_1& +& B_1\\
			\vdots &  & \ddots& &\vdots& & \vdots\\
			(a_{m1}+A_{m1}) x_1 &  +&\cdots&+& (a_{mn}+A_{mn})x_n&\subseteq & b_m& +& B_m
		\end{array}\right. .$$
		This is equivalent with $ P x\in \B $ and $A_{ij}x_j\subseteq B_i$ for all $1\leq i\leq m$ and $1\leq j\leq n$; the latter is equivalent to 
		\begin{equation}\label{constrba}
			x_j\in B_i:A_{ij},
		\end{equation} 
		for $1\leq i\leq m,1\leq j\leq n$. For $1\leq j\leq n$ the restriction \eqref{constrba} is equivalent to \begin{equation}\label{constrf}
			x_j\in \min_{1\leq i\leq m}B_i:A_{ij} =F_j.
		\end{equation} 
		Let $ j_{1}<\dots< j_{k}$ be such that  $ \F^{(c)}= (F_{j_{1}},\dots, F_{j_{k}})^{T} $ is the constraint, with  $ k $ the constraint dimension, while $ F_{j} =\R$ for $ j\in\{1,\dots, n\}\diagdown \{j_{1},\dots, j_{k}\} $. Then \eqref{constrf} amounts to
		\begin{equation}\label{constrfjh}
			x_{j_{h}} \in F_{j_{h}}\subset \R
		\end{equation}	 
		for $ 1\leq h\leq k $. 
		We may write \eqref{constrfjh} in the form $ Kx \in  \F^{(c)} $, with $ K $ the constraint matrix. Combining with the fact that $ P x\in \B $, we conclude that \(x\) is a solution of the system $			\left(\begin{matrix}
			P\\
			K
		\end{matrix} \left|  \begin{matrix}
			\mathcal B\\
			\mathcal F^{(c)}
		\end{matrix}\right. \right) 	$. Hence $ \A|\B $ and $			\left(  \begin{matrix}
			P\\
			K
		\end{matrix} \left|  \begin{matrix}
			\mathcal B\\
			\mathcal F^{(c)}
		\end{matrix}\right. \right)	$ are equivalent.	
	\end{proof}
	
	\subsection{Increasing row-echelon form}\label{subsecrow}
	Theorem \ref{rowex} states that every flexible system with real coefficient matrix can be put into increasing row-echelon form. However column permutations may be needed, so the variables may appear in different order.

	\begin{theorem}\label{rowex}
		Let $ P\in \M_{m,n} (\R)$ be of non-zero rank and $ \B\in \E^{m} $. Then there exists a permutation matrix $ H\in  \M_{n} (\R)$ such that the system $ Px\in \B $ is $ H $-equivalent to a system $Qy\in\mathcal{C} $ which is  in increasing row-echelon form  and obtained by Gaussian elimination, where $ Q\in \M_{m,n} (\R)$ is of the same rank as $ P $, $ \mathcal{C} \in \E^{m} $ and $ y=Hx $. 
	\end{theorem}
	
	The proof of Theorem \ref{rowex} is based on the following lemma and its generalization. They imply that in a reduced matrix the Gaussian operation of adding a multiple of one row to another does not change the set of real admissible solutions, provided on this row we can take a pivot equal to $ 1 $, and the neutrix at the right-hand side is minimal.
	\begin{lemma}\label{addrow1}
		Consider the reduced flexible system with real coefficients
		\begin{equation}\label{amhtd}
			\left\{	\begin{array}{rrrrrrlrr}
				x_1+a_{12}x_2+&\cdots&+a_{1n}x_n&\in& b_1&+&B_{1}\\
				a_{21}x_1+a_{22}x_2+& \cdots& +a_{2n}x_n&\in& b_2&+&B_{2}.
			\end{array}\right. 
		\end{equation}
		If $ B_{1}\subseteq B_{2} $, the system \eqref{amhtd} is  equivalent to the system with equal neutrices and with coefficient matrix of equal rank
		\begin{equation}
			\left\{	\begin{array}{rrrrrrrl}\label{amhtd2}
				x_1+& a_{12}x_2&+\cdots+&a_{1n}x_n&\in & b_1&+&B_{1}\\
				& (a_{22}-a_{21}a_{12})x_2&+\cdots+&(a_{2n}-a_{21}a_{1n})x_n &\in& b_2-a_{21}b_1&+&B_{2}.
			\end{array}\right.
		\end{equation}
	\end{lemma}
	\begin{proof} The Gaussian row-operation does not modify the rank of the coefficient matrix. Observe that $ |a_{21}|\leq 1 $ and $ B_{1}\subseteq B_{2} $, so
		\begin{equation}\label{Bfor}
			B_{2}\pm a_{21}B_{1}=B_{2}.
		\end{equation}
		Hence the row-operation leaves the neutrix at the right-hand side of the second row unchanged. We conclude that the neutrices at the right-hand side of the system \eqref{amhtd2} are equal to the neutrices at the right-hand side of the system \eqref{amhtd}. 
		
		To show the equivalence of the systems, assume $x=(x_1, \dots, x_n)^T$ satisfies the system \eqref{amhtd}. It follows directly from \eqref{Bfor} that $x$ satisfies the second row of \eqref{amhtd2}. Then $ x $ satisfies the system \eqref{amhtd2}, because it obviously satisfies the first row. Conversely, if $x=(x_1, \dots, x_n)^T$ satisfies \eqref{amhtd2}, again using \eqref{Bfor},	
		\begin{align*}
			&a_{21}x_1+a_{22}x_2+\cdots+a_{2n}x_n\\
			=&(a_{22}-a_{21}a_{12})x_2+\cdots+(a_{2n}-a_{21}a_{1n})x_n + a_{21}(	x_1+a_{12}x_2\cdots+a_{1n}x_n)\\		
			\in &b_2-a_{21}b_1+B_{2}+ a_{21}(b_1+B_{1})
			=b_2+B_{2}+a_{21}B_{1}=b_2+B_{2}.
		\end{align*}
		We conclude that $ x $ satisfies \eqref{amhtd}. Hence the two systems are equivalent.
		
	\end{proof}
	
	The following lemma on subtraction of rows in order to create zeros below some pivot on the principal diagonal is more general, and can be proved similarly.
	
	\begin{lemma}\label{addrow}
		Let $ P=(p_{ij})_{m\times n}\in \M_{m,n} (\R)$ be reduced and of rank $ r\geq 1 $, $ \B=(\beta_1,\dots, \beta_n)^T\in \E^{m},$ with \(\beta_i=b_i+B_i\) for \(1\leq i\leq m\). Let \(k\in \N\) be such that $ 1\leq k < r $. Assume that $ p_{ij}=0$ for  $k\leq i\leq m, 1\leq j< k $, $ p_{kk}=1 $,  and $ B_{k}\subseteq B_{i} $ for $ k\leq i \leq m$. Let $ Q=(q_{ij})_{m\times n}\in \M_{m,n} (\R)$  be defined by
		
		\begin{equation*}
			q_{ij}=\left\lbrace \begin{matrix}
				p_{ij}& 1\leq i\leq k, 1\leq j\leq n\\
				p_{ij}-p_{ik}p_{kj}& k+1\leq i\leq m, 1\leq j\leq n
			\end{matrix} \right.,
		\end{equation*}
		and the external vector $ \C =(\gamma_{1},\dots ,\gamma_{m})^T$ by 
		\begin{equation*}
			\gamma_{i}\equiv c_{i}+C_{i}=\left\lbrace \begin{matrix}
				\beta_{i}& 1\leq i\leq k\\
				\beta_{i}-p_{ik}\beta_{k}& k+1\leq i\leq m
			\end{matrix} \right..
		\end{equation*}Then $ C_{i}= B_{i} $ for $ 1\leq i\leq m $, $ q_{ij}=0$ for  $k+1\leq i\leq m, 1\leq j\leq k $, $ Q $ has the same rank as $ P $, and the systems $ P|\B $ and $ Q|\C $ are equivalent.
	\end{lemma}

	\begin{proof}[Proof of Theorem \ref{rowex}]
		Let $P=(a_{ij})_{m\times n}$ and $ r={\rm{r}}(P),r\geq 1 $. We prove the theorem by External induction, increasing stepwise the part of the system $P|\mathcal{B} $ having the desired form. We push all rows such that   $ P $ has at least one non-zero element to the upper side, and all zero rows of $ P|\B $ to below, 
		so between them there are possibly rows with zero elements within $ P $ and a non-zero right-hand element. To avoid notational complexity, we suppose that the rows with index $ 1,\dots, k  $  of $ P $ all have a non-zero element, and, if $ k<m $, the rows with index $ k+1,\dots, l  $ have a non-zero right-hand element, with $ l\leq m $. 
		We thus obtain a system $ P^{(0)}x\in \B^{(0)} $ which is equivalent with $ P|\mathcal{B}  $, with the same rank for the coefficient matrix.
		
		For $ 1\leq i\leq k $, we let $ |\oa_i|=\max_{1\leq j\leq n}|a_{ij}| $, choose $ a_{i\overline{j}} $  such that $ |a_{i\overline{j}}|=|\oa_i |$, and divide row $ i $ of $ P^{(0)}| \B^{(0)} $ by $ a_{i\overline{j}} $. Among the first $ k $ rows we choose a row such that the neutrix part at the right-hand side is minimal
		and permutate it with the first row. Some coefficient on the new first row is equal to $ 1 $, and we permutate the corresponding column with the first column. Then we apply Gaussian elimination to the part below the new first element.  The resulting system will be denoted by $ P^{(1)}x^{(1)}\in\B^{(1)} $, where $ x^{(1)}=H^{(1)}x$, with $ H^{(1)}\in \M_{n}(\R) $ a permutation matrix. It is reduced and in increasing row-echelon form as far as the first row is concerned, the first column of $ P^{(1)} $ has zero elements below the pivot, and the elements of its remaining columns are all limited. It follows from Lemma \ref{addrow} that the Gaussian operations used lead to an equivalent system with equal rank for the coefficient matrix and equal neutrices at the right-hand side. Hence $ \B^{(1)}_{1}\subseteq\B^{(1)}_{i}$ for $ 2\leq i \leq k$, $ \rr\left(P^{(1)}\right)=\rr(P) $ and $ P^{(1)}|\B^{(1)} $ and $ P|\B $ are $ H^{(1)} $-equivalent.
		
		Suppose that $ s<r $ and the system $P^{(s)}x^{(s)}\in\B^{(s)} $ is reduced and in increasing row-echelon form up to row $ s$, where $ x^{(s)}=H^{(s)}x$, with $ H^{(s)}\in \M_{n}(\R) $ a permutation matrix, and such that the system is $ H^{(s)} $-equivalent to $ P|\B $, with  $ \rr\left(P^{(s)}\right)=\rr(P) $, while the elements in the first $ s $ columns of $P^{(s)}$ below row $ s $ are zero, and its remaining elements are limited. 
		We insert the rows of $P^{(s)}|\B^{(s)} $ such that $ P^{(s)} $ has zero coefficients into the existing group of rows with non-zero neutrix at the right-hand side. Then the rows of  $ P^{(s)} $ are non-zero up to $ k^{(s)} $, say. We repeat the procedure sketched above, and start by constructing a reduced coefficient matrix by dividing the rows of $P^{(s)}|\B^{(s)} $ below row $ s $ by an element which in absolute value is maximal; note that this element is at most limited, so its inverse  is not an absorber of the neutrices at the right-hand side. Hence the neutrices of the right-hand side up to $ s $  continue to be contained in the neutrices of the following rows; one of these rows with minimal neutrix will be interchanged with row $ s+1 $. Then we permute columns, such that the first element of row $ s+1 $ equal to $ 1 $ occurs in column $ s+1 $. We apply Gaussian elimination to the column $ s+1 $ below row $ s+1 $. The resulting system will be denoted by $ P^{(s+1)}x^{(s+1)}\in\B^{(s+1)} $, where $ x^{(s+1)}=Sx^{(s)}$ for some  permutation matrix $ S\in \M_{n}(\R) $. The column of $ P^{(s+1)} $ below row $ s+1 $ has only zero elements. The submatrix of $ P^{(s+1)} $ below and to the right of this element is limited. It follows from the induction hypothesis and Lemma \ref{addrow} that $ \B^{(s+1)}_{i}=\B^{(s)}_{i}\subseteq\B^{(s)}_{i+1}=\B^{(s+1)}_{i+1}\subseteq\B^{(s+1)}_{s+1}$ for $ 1\leq i < s$. Hence $ P^{(s+1)}|\B^{(s+1)} $ is in increasing row-echelon form up to the $ s+1^{th} $ row. By construction and again by Lemma \ref{addrow} we have $ \B^{(s+1)}_{s+1}\subseteq\B^{(s+1)}_{i}$ for $ s+2\leq i \leq k^{(s)}$, $ \rr(P^{(s+1)})=\rr(P^{(s)}) $ and  $ P^{(s+1)}x^{(s+1)}\in\B^{(s+1)} $ is $ S $-equivalent to $ P^{(s)}x^{(s)}\in\B^{(s)}$. This implies that $ \rr(P^{(s+1)})=\rr(P)$ and that $ P^{(s+1)}x^{(s+1)}\in\B^{(s+1)} $ is $ H^{(s+1)} $-equivalent to $ Px\in\B $, with the permutation matrix $ H^{(s+1)}\equiv S H^{(s)} $. 
		
		By External induction we may thus continue up to row $ r $, and obtain a system $Qy\in\mathcal{C} $ in increasing row-echelon form up to row $ r $, where  $ y=Hx $ for some permutation matrix $ H\in \M_{n}(\R)$, such that $Qy\in\mathcal{C} $ is $ H $-equivalent to $ Px\in\B $, with $ Q=(q_{ij})_{m\times n}\in \M_{m,n} (\R)$ of the same rank as $ P $ and $ \mathcal{C} \in \E^{m} $. Observe that the elements $ q_{ij}$ for  $i>r $ and $ j<r $ are zero by the induction hypothesis, the elements below $ q_{rr} $ are zero by construction, and then the elements $ q_{ij}$ for  $i>r $ and $ j>r $ must be also zero, otherwise $ \rr(Q)>r $, a contradiction. Hence the system $Q|\mathcal{C} $ is in increasing row-echelon form.
	\end{proof}

	\subsection{On consistency}\label{subseccons}
	In this subsection we give first a criterion for consistency for a system with real coefficient matrix in row-echelon form. We apply it to obtain a critierion for an arbitrary flexible system.
	
	\begin{proposition}\label{suffreal} Let $ Q\in \M_{m,n}(\R) $ be in row-echelon form and $ \C=c+C=\linebreak(\gamma_{1},\dots,\gamma_{m})^T\in \E^{m} $.  Assume $ {{\rm{r}}}(Q)=r $. Then $ Q|\C $ is consistent if and only if $ \gamma_{j} $ is neutricial for $ r+1\leq j\leq m$, and then the systems $ Q|\C $ and $ Q^{rn}|\C^{r} $ are equivalent.
		
	\end{proposition}
	
	\begin{proof}
		Assume that  $ Q|\C $ is consistent. If $ \gamma_{j} $ is zeroless for some $ j $ with $ r+1\leq j\leq m$, we would have $ 0\in \gamma_{j} $ at row $ j $, a contradiction. Hence  $\gamma_j\equiv C_{j} $ is neutricial. Then for $ r+1\leq j\leq m$ all rows are of the form $ 0\in C_{j} $, which is automatically satisfied. Hence the systems \(Q|\C\) and $ Q^{rn}|\C ^{r} $ are equivalent.
		
		Conversely,  if $ \gamma_{j}\equiv C_{j} $ is neutricial for $ r+1\leq j\leq m$, the corresponding rows are all of the form $ 0\in C_{j} $, which as we saw is always satisfied. Hence the system \(Q|\C\) is equivalent to the remaining system  $ Q^{rn}|\C ^{r} $. Because ${\rm{r}}\left( Q^{rn}\right)  =r $, the system \(Q^{rn}|c^{r}\) is consistent, hence also $ Q^{rn}|\C^{r} $. By equivalence, \(Q|\C\) is consistent. 
	\end{proof}
	
	We now characterize the consistency of the original flexible system $ \A|\B $.

	\begin{theorem}\label{theoremsingular2}
		Let $ \A\in \M_{m,n} (\E)$, $ \B\in \E^{m} $ and $ \A|\B $ be a flexible system. Let $ P\in \M_{m,n} (\R)$ be a representative matrix of $ \A $.  Let $ \F^{(c)} $ be the constraint, $ k\leq n $ be the constraint dimension and $ K\in \M_{k,n} (\R)$ be the constraint matrix. Let ${\rm{r}}\left(  \begin{matrix}
			P\\
			K
		\end{matrix}\right)=r$. Assume that  $ Q\in \M_{m+k,n} (\R)$,  $ \C\in\E^{m+k} $ and $ Q|\C $ is a system in  increasing row-echelon form associated to $ \A|\B $ by $ P $. Then
		
		\begin{enumerate}
			\item \label{PKQ} There exists a permutation matrix $ H\in \M_{n}(\R) $ such that the systems $ \A|\B $ and $ Q|\C $ are $ H $-equivalent.
			\item \label{rQr} $ {\rm{r}(Q)}=r $.
			\item \label{rNr} The system $ \A|\B $ is consistent if and only if  $ \C^{m+k-r} $ is neutricial.
		\end{enumerate}
	\end{theorem}
	
	\begin{proof}
		\begin{enumerate}
			\item 		By Theorem \ref{thdec} the system $ \A|\B $ is equivalent to  $\left( \begin{matrix}
				P\\
				K 
			\end{matrix}\left|   \begin{matrix}
				\mathcal B\\
				\F^{(c)}
			\end{matrix}\right.\right)$.
			By Theorem \ref{rowex}   there exists a permutation matrix $ H\in \M_{n}(\R) $ such that the systems $\left( \begin{matrix}
				P\\
				K
			\end{matrix}\left|   \begin{matrix}
				\mathcal B\\
				\F^{(c)}
			\end{matrix}\right.\right)  $ and $ Q|\C $ are $ H $-equivalent. Hence $ \A|\B $ is $ H $-equivalent to $ Q|\C $. 
			\item Also by Theorems \ref{thdec} and \ref{rowex} we have $ {\rm{r}(Q)}={\rm{r}}\left( \begin{matrix}
				P\\
				K
			\end{matrix}\right)=r$. 
			\item By Part \ref{rQr} it holds that $ {\rm{r}(Q)}=r$. Then by  Proposition \ref{suffreal} the system $ Q|\C $ is consistent if and only if  $ \C^{m+k-r} $ is neutricial. This implies Part \ref{rNr}.
		\end{enumerate}
	\end{proof}

	\begin{example}\label{exinc} Consider the flexible system 
		\begin{equation}		
			\label{vd2}\left\{\begin{array}{rlllll}
				(-1+\eps\oslash)x_1+(1+\oslash)x_2 +(-0.2+\eps\pounds)x_3&\subseteq  2+\eps\pounds\\
				(1+\eps\pounds)x_1+(-1+\eps\pounds)x_2+(0.1+\oslash)x_3&\subseteq 1+\eps\oslash\\
				(1+\oslash)x_1+(-1+\oslash)x_2+(0.15+\eps\oslash)x_3&\subseteq -0.5+\oslash 
			\end{array}\right. 
		\end{equation}
		The coefficient matrix of \eqref{vd2} has the  same representative matrix as \eqref{exsys1}, having bigger neutrices. With this representative matrix the integrated system becomes
		$$\left( \begin{array}{rrrlllll}
			-1& 1&-0.2&| 2+\eps\pounds\\
			1&-1&0.1&|1+\eps\oslash\\
			1&-1&0.15&|-0.5+\oslash\\
			1&0&0 &| \oslash\\
			0&1&0&|\eps\pounds\\
			0&0&1&|\eps\pounds 
		\end{array}\right) .$$ 
		Putting the second, fifth and sixth row on top, and applying Gaussian elimination we find the system in increasing row-echelon form
		$$
		\left(\begin{array}{rrrlllll}
			1&-1&0.1&|1+\eps\oslash\\
			0&1&0&|\eps\pounds\\
			0&0&1&|\eps\pounds \\
			0& 0&0&| 3+\eps\pounds\\	
			0&0&0&|1.5+\oslash \\
			0&0&0&|-1+\oslash
		\end{array}\right) 
		$$
		The rank of the coefficient matrix is \(3\). We see that the fourth, fifth and sixth component of the right hand side are zeroless, so by Theorem  \ref{theoremsingular2} the system is inconsistent. 
	\end{example}
	\subsection{Extended parameter method}\label{subsecpar}
	
	We will now solve the system \eqref{realext}. A system with rank equal to the number of equations is solved by the parameter method  in \cite{Namopt}, which admits a solution in closed form. In the case of a system $ P|\B $, where $ P\in \M_{m,n}(\R) $ is a real coefficient matrix and $ \B=b+B $ is an external vector, the parameter method is as follows. Let $ B=(B_{1},\dots,B_{m})^{T} $. We let $ s=(s_{1},\dots,s_{m})^T $, where $ s_{i} $ is a real parameter such that $ s_{i}\in \B_{i} $ for $ 1\leq i\leq m $. We solve $ P|(b + s) $ with common methods of linear algebra, and in the end we substitute the $  s_{i}$ by their range $ B_{i} $.
	
	We will see that the parameter method also works for a system in increasing row-echelon form $ Q|\C $. So the system \eqref{realext} can also be solved, after the transformation into an equivalent system in increasing row-echelon form, as described in the proof of Theorem \ref{rowex}.
	
	Next theorem presents the solution in closed form in the case that $ P^{(m)} $ is non-singular. In addition to representatives of bounded scalar neutrices we have also parameters ranging over linear spaces of one dimension. The solution is exact.
	\begin{theorem}\label{thesolhom}
		Let $ 1\leq m\leq n\in \N $, $ P=(a_{ij})_{m\times n}\in \M_{m,n}(\R) $ be of rank $ m $ and such that $ P^{(m)} $ is non-singular, and $ \B=b+B\in \E^{m} $. Let $ \xi $ be the solution of the system $ P|\B $. Let 
		
		\begin{equation}\label{SV}
			V=\sum\limits_{k=m+1}^n \R\begin{pmatrix} -\left(P^{(m)}\right)^{-1} a_{k}^T \\
				e^{n-m}_{k-m}\end{pmatrix},
		\end{equation}  
		
		\begin{equation}\label{SW}
			W=\sum\limits_{i=1}^m \begin{pmatrix} B_i\left(P^{(m)}\right)^{-1}e_i^{m}\\
				0\end{pmatrix}.
		\end{equation}
		Then \begin{enumerate}
			\item \label{thesolhom1}$ V $ is the  solution set of the system \(P|0\).
			\item  \label{thesolhom2} \(P W=B\).
			\item  \label{thesolhom3} $ W+V $ is the  solution set of the system \(P|B\).
			\item  \label{thesolhom4}\(\xi=\left(\begin{array}{c}
				
				\left( P^{(m)}\right)^{-1}b\\ 0\end{array} \right)+W+V\).
			
			\item \label{thesolhomlm5} The solution $ \xi $ is exact.
		\end{enumerate}
	\end{theorem}
	
	\begin{proof} The set $ \xi $ is non-empty, because  $ P|b $ is consistent. 
		\begin{enumerate}
			\item It is obvious from linear algebra that $ V $ is the  solution  set of \(P|0\). 
			
			\item  Let $ s_{i}\in \R $ for $ 1\leq i\leq m $,  $ s\equiv s_{1}e^{m}_{1}+\cdots +s_{1}e^{m}_{m} $ and  
			\begin{equation*}\label{forw}
				w\equiv \begin{pmatrix}
					\left(P^{(m)}\right)^{-1} s\\0
				\end{pmatrix}\ =\sum\limits_{i=1}^m \begin{pmatrix} s_i(P^{(m)})^{-1}e_i^{m}\\0\end{pmatrix}. 
			\end{equation*}
			Assume first that $ s\in B $. Then  $ s_{i}\in B_{i} $ for $ 1\leq i\leq m $, hence $ w\in W $, while $ Pw=s $. Hence $ PW\supseteq B $. Conversely, 	let \(s\in PW\). Then there exists \(w\in W\) such that 	\(Pw=s\). Because \(w\in W\), by \eqref{SW} it holds that \(s\in B\). Hence \(PW\subseteq B\). We conclude that $ PW= B $.
			\item	Let $ \xi_{B} $ be the solution set of $ P|B $. Clearly $ P(W+V) =PW+PV=B+0=B$, so $ W+V\subseteq  \xi_{B}$. On the other hand, let $ s\in B $ and $ \xi_{s} $ be the solution of $ P|s $. By Part \ref{thesolhom2} there exists $ w_{s}\in W $ such that $ Pw_s=s $. 
			Then it follows from  linear algebra that $ \xi_{s}=w_s+V $ is the solution of $ P|s $. Hence $ \xi_{B}=\cup_{s\in B}\xi_{s}=\cup_{s\in B} (w_{s}+V)\subseteq W+V $. We conclude that $ \xi_{B}= W+V $.
			\item Clearly  $ \left( \begin{array}{c}		 	
				\left( P^{(m)}\right)^{-1}b\\ 0\end{array} \right) \in \xi$. Then Part \ref{thesolhom4} follows from Theorem \ref{theoremsol} and Part \ref{thesolhom3}.
			\item The matrix $ P $ has real coefficients. Then by distributivity, Part \ref{thesolhom4}, Part \ref{thesolhom2} and Part \ref{thesolhom1} it holds that 
			\begin{align*}
				P\xi&=P\left( \left(\begin{array}{c}				
					\left( P^{(m)}\right)^{-1}(b)\\ 0\end{array} \right)+W+V\right) \\
				=&P \left( \begin{array}{c}				
					\left( P^{(m)}\right)^{-1}(b)\\ 0\end{array} \right) +PW+PV			=b+B+0=b+B.
			\end{align*} 
			Hence the solution $ \xi $ is exact.
		\end{enumerate}
	\end{proof}

	\section{The neutrix part of the solution of a system with real coefficient matrix}\label{secprop}

	Let $ P|\B $ be a flexible system, where $ P\in \M_{m,n}(\R) $ is a real coefficient matrix of rank $ m $. We will see that the solution set $ \xi $ is exact. Its neutrix $ {\rm{N}}(\xi) $ has the form of a direct sum of its linear part $ \xi_{(L)} $	and a modular part $ \xi_{(M)} $. The dimension of $ \xi_{(L)} $ depends on the rank of $ P $ and the dimension of $ \xi_{(M)} $ on the number of non-zero neutrices at the right-hand side.
	
	\begin{theorem}\label{thesolhomlm}	
		Let $ 1\leq m\leq n\in \N $, $ P=(a_{ij})_{m\times n}\in \M_{m,n}(\R) $ be of rank $ m $ and such that $ P^{(m)} $ is non-singular, and $ \B=b+B\in \E^{m} $. Let $ \xi $ be the solution of the system $ P|\B $. Let $ V $ be given by 	\eqref{SV} and $ W $  by \eqref{SW}.			Then
		\begin{enumerate}
			\item \label{thesolhomlm1}$ V $ is a linear space, equal to the direct sum 
			\begin{equation*}\label{solrealparmn5}
				V=\oplus_{m+1\leq k \leq n}  \R\begin{pmatrix} -\left(P^{(m)}\right)^{-1} a_{k}^T\\
					e^{n-m}_{k-m}\end{pmatrix},			
			\end{equation*}
			with dimension \begin{equation*}\label{solrealparmnl5}
				\dim V=n-m.
			\end{equation*}
			\item  \label{thesolhomlm2} \(W\) is a bounded neutrix, equal to the direct sum 
			\begin{equation}\label{solrealparmnl3}
				W= \oplus_{B_k\supset\{0\}} B_k\begin{pmatrix} \left(P^{(m)}\right)^{-1}e_k^{m}\\
					0\end{pmatrix},
			\end{equation} with dimension
			\begin{equation}\label{soldimBl}
				\dim(W)=\sharp\{k|B_k\supset\{0\}\}.	
			\end{equation} 			
			\item  \label{thesolhomlm3}  \({\rm{N}}(\xi)=W\oplus V\).			
			
			\item \label{thesolhomlm4} $ V $ is the linear part of $ \xi $, and $ W $ is a modular part of $ \xi $.

		\end{enumerate}
	\end{theorem}
	
	\begin{proof} 
		\begin{enumerate}
			\item Being the solution of $ P|0 $, the set $ V $ is a linear space. The set of vectors
			\begin{equation*}		
				G_{V}\equiv	\left\lbrace  \begin{matrix} -\left(P^{(m)}\right)^{-1} a_{m+1}^T\\
					e^{n-m}_{1}\end{matrix},\dots,\begin{matrix} -\left(P^{(m)}\right)^{-1} a_{n}^T\\
					e^{n-m}_{n-m}\end{matrix} \right\rbrace  
			\end{equation*}
			is linearly independent. Hence
			\begin{equation*}
				V=\sum\limits_{k=m+1}^n \R\begin{pmatrix} -\left(P^{(m)}\right)^{-1} a_{k}^T\\
					e^{n-m}_{k-m}\end{pmatrix}=\oplus_{m+1\leq k \leq n} \R\begin{pmatrix} -\left(P^{(m)}\right)^{-1} a_{k}^T\\
					e^{n-m}_{k-m}\end{pmatrix}
			\end{equation*} and \(\dim V=n-m\). 
			\item 
			
			Formula \eqref{solrealparmnl3} is a consequence of the fact that the set of $m $ vectors	
			\begin{equation*}		
				G_{W}\equiv		\left\lbrace  \begin{matrix} -\left(P^{(m)}\right)^{-1} e_{1}^{m}\\
					0\end{matrix},\dots,\begin{matrix} -\left(P^{(m)}\right)^{-1} e_{m}^{m}\\
					0\end{matrix} \right\rbrace  
			\end{equation*}
			is linearly independent. Being a direct sum of scalar neutrices, the set $ W $ is a neutrix. By Remark \ref{remboun} the components of the neutrix $ B $ are bounded. Then it follows from \eqref{solrealparmnl3} that the components of the neutrix $ W $ are bounded, hence also $ W$ is bounded. Formula \eqref{soldimBl} follows from the fact that $ B_{k}\left(P^{(m)}\right)^{-1} e_{k}^{m}\supset\{0\} $ if and only if $ B_{k}\supset\{0\} $, for $ 1\leq k\leq m $.

			\item  	Clearly \(G_{V}\cup G_{W}\) is linearly independent. This implies that \({\rm{N}}(\xi)=V \oplus W\). 
			
			\item By Part \ref{thesolhomlm3} it holds that \({\rm{N}}(\xi)=V \oplus W\), while $ V $ is a linear space by Part \ref{thesolhomlm1}, and $ W $ is a bounded neutrix by Part \ref{thesolhomlm2}. Then Theorem \ref{decNLnew} implies that $ V $ is the linear part of $ {\rm{N}}(\xi) $. Then it follows from Definition \ref{deflinpar} that $ V= \xi_{(L)}$ and  that $  W$  is a modular part of $\xi$.

		\end{enumerate}
	\end{proof}

	\section{Invariance of rank of integrated matrices}\label{secrank}
	
	Representative matrices of an external matrix may have different ranks; this is obvious for a neutricial matrix, which has both a zero representative matrix and non-zero representative matrices. In contrast, Proposition \ref{Namproof} shows that the rank of the coefficient matrix of an integrated system associated to a flexible system is always the same.
	
	\begin{proposition}\label{Namproof} Consider the system $ \A|\B $, where $ \mathcal{A} \in \M_{m,n}(\E)$ and $ \B\in \mathcal{\E}^{m}$ with $ {\rm{N}}(\B)=B $. Let $			\left( \begin{matrix}
				P\\
				K
			\end{matrix} \left|  \begin{matrix}
				\mathcal B\\
				\mathcal F^{(c)}
			\end{matrix}\right. \right)	$ and $			\left( \begin{matrix}
			P'\\
			K
		\end{matrix} \left|  \begin{matrix}
		\mathcal B\\
		\mathcal F^{(c)}
	\end{matrix}\right. \right)	$ be two associated integrated systems, where \(P, P'\) are two representative matrices of \(\mathcal{A} \) and $ K $ is the constraint matrix. Let $ r={\rm{r}}\left(  \begin{matrix}
			P\\
			K
		\end{matrix}  \right) $ and $ r'={\rm{r}}\left(  \begin{matrix}
			P'\\
			K
		\end{matrix}\right) $. 		
		Let $ C\in {\mathcal{N}}^{m}, C'\in {\mathcal{N}}^{m}$, and $ Qy\in C $  be a system in  increasing row-echelon form associated to $ \mathcal{A}x\subseteq B $ by $ P $ and $ Q'z\in C' $  be a system in increasing row-echelon form associated to $ \mathcal{A}x\subseteq B $ by $ P' $, where $ y=Hx $ and $ z=H'x $ for some permutation matrices $ H,H'\in \M_{n}(\R) $.  Let \(\xi\) be the solution of the homogeneous system $
		\left( \begin{matrix}
			P\\
			K
		\end{matrix} \right) x \in \left( \begin{matrix}
			B\\
			\mathcal F^{(c)}
		\end{matrix}\right) 
		$ and $ \xi'$ be the solution  of the homogeneous sytem $
		\left( \begin{matrix}
			P'\\
			K
		\end{matrix}\right) x\in \left( \begin{matrix}
			B\\
			\mathcal F^{(c)}
		\end{matrix} \right)
		$. 
\\	
		Then 
		\begin{enumerate}
			\item \label{forlins} \(\xi_{(L)}=\xi'_{(L)}\)
			\item \label{forsols} $ H\xi $ is the solution of \(Q^{rn}|C^{r}\) and $ H'\xi' $ is the solution of \(\left(Q'\right)^{r'n}|\left(C'\right)^{r'}\).
			
			\item \label{forranks} $ 
			r	={\rm{r}}(Q)={\rm{r}} (Q')=r'$.
			
		\end{enumerate}
	\end{proposition}
	
	\begin{proof} 
		Being homogeneous, the system $ \mathcal{A}x\subseteq B $ is consistent. Then its solution is a neutrix $N$.  By Theorem \ref{thdec} the homogeneous systems $
		\left( \begin{matrix}
			P\\
			K
		\end{matrix} \right) x \in \left(    \begin{matrix}
			B\\
			\mathcal F^{(c)}
		\end{matrix}\right)
		$ and $
		\left( \begin{matrix}
			P'\\
			K
		\end{matrix}\right) x\in \left(    \begin{matrix}
			B\\
			\mathcal F^{(c)}
		\end{matrix} \right)
		$ are both equivalent to $ \mathcal{A}|\B $, hence they are equivalent. Consequently,
		\begin{equation}\label{solH}
			\xi=N=\xi'.
		\end{equation}
		\begin{enumerate}
			\item From \eqref{solH} we derive that \(\xi_{(L)}=\xi'_{(L)}\).
			
			\item 	It follows from Theorem \ref{rowex} that $ H\xi $ is the solution of $ Qy\in C $ and $ H'\xi' $ is the solution of $ Q'z\in C' $. Again by Theorem \ref{rowex}
			\begin{equation}\label{forrr}	
				{\rm{r}}(Q)={\rm{r}}\left( \begin{matrix}
					P\\
					K
				\end{matrix}\right)=r , \quad {\rm{r}}(Q')={\rm{r}}\left( \begin{matrix}
					P'\\
					K
				\end{matrix}\right) =r'.
			\end{equation}
			Then it follows from Proposition \ref{suffreal} that $ H\xi $ is the solution of \(Q^{rn}|C^{r}\) and $H' \xi' $ is the solution of \(\left(Q'\right)^{r'n}|\left(C'\right)^{r'}\). 
			
			\item 	By Parts \ref{forlins} and \ref{forsols}
			\begin{align}\label{forrqr}
				{\rm{r}}(Q)&={\rm{r}}\left(Q^{rn}\right) =n-\dim\left((H\xi)_{(L)}\right)=n-\dim\left(H\xi_{(L)}\right)=n-\dim\left(\xi_{(L)}\right)\notag\\
				&=n-\dim\left(\xi'_{(L)}\right)=n-\dim\left(H'\xi'_{(L)}\right)=n-\dim\left((H'\xi')_{(L)}\right)\notag\\
				&={	\rm{r}}\left(\left(Q'\right)^{r'n}\right)
				={\rm{r}}\left(Q'\right).
			\end{align}

			Formulas \eqref{forrr} and \eqref{forrqr} imply Part \ref{forranks}. 
		\end{enumerate}
	\end{proof}
	\section{Proof of the Main Theorem}\label{secproof}
	\begin{proof}[Proof of Theorem \ref{maintheorem}]
		
		\ref{part1}. By Theorem \ref{thdec} the system \(\A x\subseteq \B\) is equivalent with		$ \left(  \begin{matrix}
			P\\
			K
		\end{matrix}\right) x\in \left(   \begin{matrix}
			\mathcal B\\
			\mathcal F^{(c)}
		\end{matrix}\right)$, 
		where $ \left( \begin{matrix}
			P\\
			K
		\end{matrix}\right) \in \M_{m+k,n}(\R) $ and $\left( \begin{matrix}
			\mathcal B\\
			\mathcal F^{(c)}
		\end{matrix}\right) \in \E^{{m+k}}$. 		
		By Theorem \ref{rowex} there exists a permutation matrix  \(H\in \M_{n}(\R)\) such that  the system $ \left( \begin{matrix}
			P\\
			K
		\end{matrix}\right) x\subseteq \left( \begin{matrix}
			\mathcal B\\
			\mathcal F^{(c)}
		\end{matrix}\right)$ 				
		is \(H\)-equivalent with a system \(Q y\in \mathcal{C}\) which is  in increasing row-echelon form and obtained by Gaussian elimination, where $ Q\in \M_{m,n} (\R)$, $ {\rm{r}} (Q)=r $ and $ \mathcal{C}\equiv c+C \equiv (\gamma_{1},\dots,\gamma_{m+k})^{T} \in \E^{m+k}$, with $ \gamma_{i}=c_{i}+C_{i} $ for $ 1\leq i\leq m +k$. Hence also
		\(\A x\subseteq \B\) is $ H $-equivalent with the system \(Qy\subseteq \C\).

		\ref{part2}. Since $ {\rm{r}} (Q)=r $, it follows from Proposition \ref{suffreal}
		that the system $ Q|\C $ is consistent if and only if $ \gamma_{j} $ is neutricial for $ r+1\leq j\leq m+k$.
		
		\ref{part3}.   By Proposition \ref{suffreal} the consistent system $ Q|\C $ is equivalent to the system \(Q^{rn}|\C^{r}\). Because $ Q^{(r)} $ is an upper triangular matrix of dimension $ r\times r $, so is $ \left(Q^{(r)}\right)^{-1}$. Put 
		\begin{equation*}\label{QV}
			W=\sum\limits_{i=1}^r \begin{pmatrix} C_i\left(Q^{(r)}\right)^{-1}e_i^{r}\\
				0\end{pmatrix}
		\end{equation*}
		\begin{equation*}\label{QW}
			V=\sum\limits_{k=r+1}^n \R\begin{pmatrix} -\left(Q^{(r)}\right)^{-1} (q^{r}_{k})^T\\
				e^{n-r}_{k-r}\end{pmatrix}
		\end{equation*}
		
		Then \eqref{parasingular1} follows from Theorem \ref{thesolhom}.\ref{thesolhom4}. By Part \ref{part1} the solution of $ \A|\B $ is given by $ \xi=H^{-1} \zeta$. 
		
		\ref{part4}. By Theorem \ref{thesolhomlm}.\ref{thesolhomlm4} it holds that $ \zeta_{(L)}=V $, and $ \dim(\zeta_{(L)}) =m+k-r$ by Theorem \ref{thesolhomlm}.\ref{thesolhomlm1}.
		
		\ref{part5}. By Theorem \ref{thesolhomlm}.\ref{thesolhomlm2} it holds that $ \zeta_{(M)}\equiv W $ is a bounded neutrix and $ \dim(\zeta_{(M)})=\sharp\{i\leq r|C_i\supset\{0\}\}$.
		
		\ref{part6}  
		By Proposition \ref{Namproof}.\ref{forranks} the rank of  a coefficient matrix of an associated integrated system does not depend on the choice of a representative matrix. Let the neutrix $ N $ be the solution of the homogeneous system $ \A x\subseteq B  $, with \(B\) the neutrical vector associated to the external vector \(\B\). Then $ {\rm{N}}(\xi)= N$ by Theorem \ref{theoremsol}, so the neutrix part of $ \xi $ does not depend on the choice of a representative matrix. Also $ \xi_{(L)}=N_{(L)} $, so the linear part of $ \xi $ does not depend on the choice of a representative matrix.  Let  $ M $ be a modular part of $\xi$.  Then $ M $ is a modular part of ${\rm{N}}(\xi)=N$, so $  N=N_{(L)}\oplus M $. By Theorem \ref{themod}.\ref{themod2} we have $ \dim(M )=\dim(N) -\dim(N_{(L)})$, which does not depend on the choice of a representative matrix.
	\end{proof}
	\section{Essential parts and feasible systems}\label{secneg}
	
	Consider a flexible system $ \A|\B $. The appearence of neutrices in the coefficient matrix induces feasibility equations in the associated integrated system $ \left(  \begin{matrix}
		P\\
		K
	\end{matrix} \left|  \begin{matrix}
		\mathcal B\\
		\mathcal F^{(c)}
	\end{matrix}\right. \right)  $.  Sometimes they just restrict the range of some the variables, such as in Example \ref{exfeas1}, but it is also possible that they interfere with the original system, which was the case in Example \ref{exinc}. Indeed, the solution strategy of the Main Theorem may involve change of rows, so a row of the "constraint part" $ K|\F^{(c)} $ could be inserted into the "representative part" $ P|\B $. We call a system \emph{feasible} if this does not need to happen, otherwise said, if $ \A|\B $ is equivalent with  $P|\B $. The equivalence will  be a consequence of a more general property, which divides a system $ \A|\B $ into an "essential part" and a "remaining part" of inclusions which may be neglected. Proposition \ref{suffreal} is a special case of this property, indicating that rows with zeros in the coefficient matrix and a neutrix at the right-hand side may be omitted.

	To start with we recall the notion of determinant and non-singularity given in \cite{Jus} and introduce some notation. 
	
	For $ \mathcal{A}\in \M_{n}(\E) $, the determinant $%
	\Delta \equiv \det(\mathcal{A}) \equiv d+D$ is defined in the usual way through sums of signed products. Also minors are defined in the usual way.

	\begin{definition}\label{defmatnonsin}
		Let $ \A\in \M_{n}(\E) $. Then $\mathcal{A}$ is called {\em non-singular} if $\Delta $ is zeroless.
	\end{definition}

	Observe that a representative matrix of a non-singular matrix $ \A $ is always non-singular.

	\begin{notation}
		Consider the flexible system $\A|\B$, where $\A\in \M_{m,n}(\E)$ and $\B\in \E^{m}$. Let $1\leq r<m$. For each $j$ such that $1\leq j\leq n$ we write $\overline{{A}^{r}_{j}}=\max\limits_{1\leq i\leq r} A_{ij}$, $\overline{{A}^{m-r}_{j}}=\max\limits_{1+r\leq i\leq m} A_{ij}$, $\overline{{B}^{r}}=\max\limits_{1\leq i\leq r} B_i$  and   $\underline{{B}^{m-r}}=\min\limits_{1+r\leq i\leq m} B_i$.
	\end{notation}

	\begin{definition}\label{defess}
		Let $ \A\in \M_{m,n} (\E)$, $ \B\in \E^{m} $ and $ \A|\B $ be a flexible system. Let $ \A_{E}| \B_{E}$ be a subsystem of $ \A|\B $.  The subsystem $ \A_{E}| \B_{E}$ is called \emph{essential} if $ \A|\B $ and $ \A_{E}| \B_{E}$ are equivalent.
	\end{definition}
	
Assume that the rows and columns of a flexible system $ \A|\B $ are ordered in such a way that the submatrix $ \A^{(r)} $ is non-singular. Theorem \ref{lemma nghiem tham so singula} below gives a criterion such that the first $ r $ rows of the system form an essential subsystem. In fact $ r $ should be the rank of a representative system and at the right-hand side the neutrices below row $ r $ should be at least as big as the neutrices up to row $ r $, while in contrast the biggest neutrix in each column of the coefficient matrix should appear above row $ r $. In addition $ \det\left(\A^{(r)} \right) $ should not be an absorber of the maximal neutrix at the right-hand up to $ r $.
	
	\begin{theorem}\label{essential} 
		 Let  $ \A=(\alpha_{ij})_{m\times n}\in \M_{m,n} (\E)$ be limited and $ \B=(\beta_1,\dots, \beta_n)^T\in \E^{m} $. Let \({\rm{N}}(\alpha_{ij})=A_{ij}\) and \(\rm{N}(\beta_i)=B_i\) for \(1\leq i\leq m, 1\leq j\leq n\). 
		Let $r\in \N$ be such that $  1\leq r<m $. Assume that:
		\begin{enumerate}
			\item $ \A^{(r)}  $ is non-singular.
			\item  \label{essetial part1} There exist a representative matrix \(P=(a_{ij})_{m\times n}\) of \(\A\) and a representative vector \(b\) of \(\B\)  such that \(\rr(P|b)=r.\)
			\item \label{gt 1 pp tham so singular1} $\det\left(\A^{(r)}\right)$ is not an absorber of $\overline{{B}^{r}}$.
			\item \label{gt 2 pp tham so singular1} $\overline{{B}^{r}}\subseteq \underline{{B}^{m-r}}$.
			\item \label{gt3} $\overline{{A}^{m-r}_{j}}\subseteq \overline{{A}^{r}_{j}}$ for all $1\leq j\leq n$. 
		\end{enumerate}
		Then $ \A^{rn}|\B^{r} $ is an essential part of $ \A|\B $.	
	\end{theorem}
	
	\begin{proof} Obviously, a vector $x=(x_1,\dots, x_n)^T\in \R^n$ which is a solution of the system $ \A|\B $ is a solution of $ \A^{rn}|\B^{r} $. 
		
		Conversely, assume that $x=(x_1,\dots, x_n)^T\in \R^n$ is a solution of  the system  $ \A^{rn}|\B^{r}  $. By assumption \eqref{essetial part1}, there exists a representative matrix $ P=(a_{ij})_{m\times n} $ of $ \A $ and a representative vector $ b=(b_1,\dots, b_n)^T $ of $ \B $ such that $(P|b)$ has rank $r$. 	Let $ r+1\leq k\leq m $. We need to prove that $x$ satisfies the $k^{th}$ inclusion  of the system $ \A|\B $, i.e.
		\begin{equation}\label{inclk}
			\sum\limits_{j=1}^n \alpha_{kj} x_j=\sum\limits_{j=1}^n (a_{kj}+A_{kj}) x_j=\sum\limits_{j=1}^n a_{kj} x_j +\sum\limits_{j=1}^nA_{kj}x_j\subseteq \beta_k.
		\end{equation}

		We prove first the neutrix part. Let $1\leq j\leq n$. Because $\overline{{A}^{m-r}_{j}} \subseteq \overline{{A}^{r}_{j}}$, $\overline{{A}^{r}_{j}}=A_{i_0j}  $ for some $i_0  $ with $1\leq i_0\leq r$ and $x$ is a solution of  $ \A^{rn}|\B^{r} $, 
		\begin{equation*}
			A_{kj} x_j\subseteq \overline{{A}^{m-r}_{j}} x_j\subseteq \overline{{A}^{r}_{j}} x_j=A_{i_0j}x_j\subseteq B_{i_0}\subseteq\overline{B^{r}}\subseteq \underline{B^{m-r}}\subseteq B_k .
		\end{equation*}
		Hence
		\begin{equation}\label{neutpart}
			\sum\limits_{j=1}^nA_{kj}x_j\subseteq B_k. 
		\end{equation}
		
		Secondly, we show that $\sum\limits_{j=1}^n a_{kj}x_j\in \beta_k$. For $ 1\leq i\leq m $ we denote the $ i^{th} $ row of $ (P|b) $ by $u_i\equiv(a_{i1}, \dots, a_{in}, b_i)$.  Because  ${{\rm{r}}}(P|b)=r$ and the matrix $ \A^{(r)} $ is non-singular, it holds that $d_r\equiv\det\left(P^{(r)}\right)\not=0$, i.e. $ P^{(r)} $ is also non-singular. Then there exist real numbers $t_1,\dots, t_r$ such that 
		\begin{align}\label{bieu dien tuyen tinh pp tham so} 
			u_k=t_1u_1+\cdots+t_ru_r;
		\end{align}	
		in fact it follows from Cramer's rule that 	$t_i=\dfrac{\det\left(P^{(i)}_{rk}\right)}{d_r}$ for $1\leq i\leq r$, where $$P^{(i)}_{rk}\equiv\begin{pmatrix}
			a_{11} & \dots & a_{1r}\\
			\vdots & \ddots &\vdots\\
			a_{(i-1)1} & \cdots & a_{(i-1)r}\\
			a_{k1}& \cdots & a_{kr}\\
			a_{(i+1)1} &\cdots & a_{(i+1)r}\\
			\vdots & \ddots & \vdots \\
			a_{r1} & \cdots & a_{rr}
		\end{pmatrix}.$$
		
		By assumption \eqref{gt 1 pp tham so singular1}, and the fact that $\det\left(P^{(i)}_{rk}\right)$ is limited, we have for $1\leq i\leq r$
		\begin{equation*}
			t_i\overline{B^{r}}=\det\left(P^{(i)}_{rk}\right)\dfrac{\overline{B^{r}}}{d_r} =\det\left(P^{(i)}_{rk}\right) \overline{B^{r}}\subseteq \overline{B^{r}}.
		\end{equation*} 
		By assumption \eqref{gt 2 pp tham so singular1}
		\begin{equation}\label{danh gia neutrix pp tham so}
			t_1B_1+\cdots +t_rB_r\subseteq t_1\overline{B^{r}}+\cdots +t_r\overline{B^{r}}\subseteq r\overline{B^{r}}=\overline{B^{r}}\subseteq\underline{B^{m-r}}\subseteq B_k.
		\end{equation}	
		Because  $x$ is a solution of the system $ \A^{rn}|\B^{r} $,  $$\left\{\begin{array}{ccccccccc}
			a_{11} x_1&+&a_{12}x_2&+&\cdots&+& a_{1n}x_n&\in & b_1+B_1\\
			\vdots & & \vdots& & \ddots& & \vdots& & \vdots\\
			a_{r1} x_1& +&a_{r2}x_2& +&\cdots&+& a_{rn}x_n&\in & b_r+B_r
		\end{array}\right.,$$
		hence
		$$\left\{\begin{array}{ccccccccc}
			t_1a_{11} x_1&+&t_1a_{12}x_2&+&\cdots&+& t_1a_{1n}x_n&\in & t_1(b_1+B_1)\\
			\vdots & & \vdots& & \ddots& & \vdots& & \vdots\\
			t_ra_{r1} x_1& +&t_ra_{r2}x_2& +&\cdots&+& t_ra_{rn}x_n&\in & t_r(b_r+B_r)
		\end{array}\right..$$	
		Consequently,
		\begin{align*}
			\big(t_1a_{11} x_1+ t_1a_{12}x_2+\cdots+t_1a_{1n}x_n\big)+\cdots +\big(t_ra_{r1} x_1+t_ra_{r2}x_2+\cdots+t_ra_{rn}x_n\big)	\in \\
			t_1(b_1+B_1)+\cdots +t_r(b_r+B_r),
		\end{align*}
		hence
		\begin{align*}
			\big(t_1a_{11}+ \cdots + t_r a_{r1}\big)x_1+ \big(t_1a_{12}+\cdots+t_ra_{r2}\big)x_2+\cdots +\big(t_1a_{1n} +\cdots+ t_ra_{rn}\big) x_n \in\\
			\big(t_1b_1+\cdots +t_rb_r\big)+\big(t_1B_1+\cdots+ t_rB_r\big).
		\end{align*}
		Then \eqref{bieu dien tuyen tinh pp tham so}  and \eqref{danh gia neutrix pp tham so} imply that 
		\begin{equation}\label{realpart}
			a_{k1}x_1+\cdots +a_{kn}x_n \in b_k+\big(t_1B_1+\cdots+ t_rB_r\big)\subseteq b_k+B_k.
		\end{equation}
		Formula \eqref{inclk} follows from \eqref{realpart} and \eqref{neutpart}. Hence $ \A^{rn}|\B^{r} $ is an essential part of $ \A|\B $.
	\end{proof}
	\begin{corollary}\label{essential2}With the notations of Theorem  
		\ref{essential}, let $ \A\in \M_{m,n} (\R)$ be limited. Assume that:
		\begin{enumerate}
			\item $ \A^{(r)}  $ is non-singular.
			\item  \label{essetial part1l} There  exists a representative vector \(b\) of \(\B\)  such that \(\rr(\A|b)=r.\)
			\item \label{gt 1 pp tham so singular1l} $\det\left(\A^{(r)}\right)$ is not an absorber of $\overline{{B}^{r}}$.
			\item \label{gt 2 pp tham so singular1l} $\overline{{B}^{r}}\subseteq \underline{{B}^{m-r}}$.
		\end{enumerate}
		Then $ \A^{rn}|\B^{r} $ is an essential part of $ \A|\B $.	
	\end{corollary}
	
	\begin{proof}
		Because \(\A\) is a real matrix, it holds that $\overline{{A}^{m-r}_{j}}= \overline{{A}^{r}_{j}}=0$ for $1\leq j\leq n$. Hence the result follows from Theorem \ref{essential}. 
	\end{proof}
	\begin{definition}\label{sysfeas}
		A flexible system $ \A|\B $ is said to be \emph{feasible} if there exists a representative matrix $ P $ of $ \A $ such that $ \A|\B $ and $ P|\B $ are equivalent.
	\end{definition}
	By Theorem \ref{thdec} a system $ \A|\B $ is equivalent with  its associated integrated system $
	\left(  \begin{matrix}
		P\\
		K
	\end{matrix} \left|  \begin{matrix}
		\mathcal B\\
		\mathcal F^{(c)}
	\end{matrix}\right. \right),
	$
	here $ P $ is a representative matrix   of $ \A $, $ \F^{(c)} $ is the constraint and $ K $ is a constraint matrix. So in case the system is feasible, the representative part $ P|\B $ is an essential part, in other words we may neglect the constraint part $ K|\F^{(c)} $. Theorem \ref{thefeas} gives conditions for this to happen. For convenience it is formulated for systems of full rank.
	\begin{theorem}\label{thefeas}
		Let $ m\leq n $, $ \A\in \M_{m,n} (\E)$ be limited, $ \B\in \E^{m} $ and $ \A|\B $ be a flexible system. 
		Let $\overline{m} $ be the maximum in absolute value of all minors of order $ m $. Assume that \(\overline{m}\neq 0\).	 Let $ \mathcal{F}^{(c)}$ be the feasibility space of $ \A|\B $ with components $F_{1},\dots,F_{k} $ and let $\underline{F} =\min\{F_{1},\dots,F_{k} \} $. If
		\begin{enumerate}
			\item \label{feas1} $\overline{{B}}\subseteq \underline{F}$
			\item \label{feas2} $\overline{m}$ is not an absorber of $\overline{{B}}$,		
		\end{enumerate} the system $ \A|\B $ is feasible.
	\end{theorem}
	
	\begin{proof}
		Because \(\A\) has a non-zero minor of order \(m\), there exists a representative matrix   $ P $ of \(\A\) and  a representative vector $ b $  of $ \B $ such that $ {\rm{r}}(P)={\rm{r}}(P|b)=m$. Let $ K $ be the constraint matrix of $ \A|\B $. By Theorem \ref{thdec} the system $ \A|\B $ is equivalent to $ \left( \begin{matrix}
			P\\
			K
		\end{matrix}\left|   \begin{matrix}
			\mathcal B\\
			\F^{(c)}
		\end{matrix}\right.\right)  $. We may change the order of appearance of the variables to obtain an $ H $-equivalent system $ \left(  \begin{matrix}
			P'\\
			K'
		\end{matrix}\left|   \begin{matrix}
			\mathcal B\\
			\F^{(c)}
		\end{matrix}\right.\right)  $ such that $ \det\left(\left(P'\right)^{(m)}\right)= \overline{m} $, where $ K' $ is again a constraint matrix and $ H \in \M_{n}(\R)$ is a permutation matrix. Then $ \det\left(\left(P'\right)^{(m)}\right) $ is  not an absorber of $\overline{{B}}$, while $\overline{{B}}\subseteq \underline{F}$. By Corollary \ref{lemma nghiem tham so singula} the system $ \left(  \begin{matrix}
			P'\\
			K'
		\end{matrix}\left|   \begin{matrix}
			\mathcal B\\
			\F^{(c)}
		\end{matrix}\right.\right)  $ is equivalent with $ P'|\B $, hence $ H $-equivalent with $ P|\B $. Hence $ \A|\B $ is feasible.
	\end{proof}
	\section{On robustness}\label{secrob}
	Informally, a property is robust if it is stable under small perturbations. Often robustness is studied in the context of optimization, and typically, if a minimum is attained at some point $ u $, one looks for a convex set $ V $ in the neighborhood of $ u $ such that for values in $ V $ the same minimum is attained, and the determination of the biggest set $ V $ becomes a maximization property, see e.g. \cite{Soyster}, \cite{ Kouvelis}, \cite{Bental2}, \cite{Bental}. Strict robustness requires that the property is totally unchanged in some neighborhood of a given value, in other cases it is only asked that the property almost holds \cite{Bertsimas}, \cite{Fischetti},  \cite{Liebchen}, and one may speak about light of recoverable robustness. In our context of the study of inclusions we choose to study a form of strict robustness, i.e., persistence of a property in a convex neighborhood.
	\begin{definition}\label{defrob}
		Let $ u\in \R^{d} $ with $ d\geq 1 $ and $ R(u) $ be a property. We say that $ R(u) $ is \emph{robust} if there exists a convex set $ V\subseteq \R^{d} $ such that $ \{u\}\subset V $ and $ R(v) $ holds for every $ v\in V $. The \emph{robustness domain} $ W $ of $ R(u) $ is defined by $ W=\{v\in \R^{d}|R(v)\} $.
	\end{definition}
	
	\begin{definition}Consider the system $ P|\B $, where $ P\in M_{m,n}(\R) $ is reduced. Let $ \mathcal{Q} \in \mathcal{M}_{m,n}(\mathcal{N})$. If the systems $ P|\B $ and $ (P+\mathcal{Q}) |\B $ are equivalent, the system $ (P+\mathcal{Q})|\B $ is called a \emph{strict perturbation} of $  P|\B$. A strict perturbation $ (P+\mathcal{Q}) |\B $ is \emph{limited}	if the matrix $ P+\mathcal{Q} $ is limited. The \emph{robustness matrix} is the maximal matrix $ \mathcal{R} $ in the sense of inclusion such that $ \mathcal{R}|\B $ is a limited strict perturbation of $  P|\B $. 
	\end{definition}
	
	Assume $   Px\in \B $, where $ P $ is a real matrix, $ x $ a real vector and $ \B $ an external vector. Let $ (P+\mathcal{Q})|\B $ be a strict perturbation of $  P|\B$, where $ \mathcal{Q}=(Q_{ij})_{m\times n}$, with $Q_{ij}  $ a neutrix for $ 1\leq i\leq m, 1\leq j\leq n $. Consider a representative matrix $ \hat{q}=(q_{ij})_{m\times n} $ of $ \mathcal{Q} $, i.e. $q_{ij}\in Q_{ij}  $ for $ 1\leq i\leq m, 1\leq j\leq n $. We may identify $ \hat{q}$ with a vector $ q\in \R^{mn} $. Consider the property \begin{equation}\label{proprub}
		R(q) := (P+\hat{q})x\in \B.
	\end{equation}
	Then  $ R(0) $ corresponds to $ Px\in \B $. We see that $ R(0) $ is robust, with $ R(q)\Leftrightarrow R(0) $ for $ q\in \oplus _{1\leq i\leq m, 1\leq j\leq n } Q_{ij}$, which is convex indeed.
	
	The neutrices of a limited strict perturbation, and in particular the  robustness matrix, must be strictly contained in $ \pounds $, so are at most equal to $ \oslash $. Note a perturbation by the neutrix $ \pounds $ tends to be too incisive, for it would lead to a coefficient matrix that cannot be put in reduced form, and in many cases to inconsistency. 
	\begin{example}
		Consider the system $ P|\B $ given by $
		P|\B=
		\left(  \begin{array}{ll|l}
			1&1&3+\oslash\\
			1&-1&1+\oslash\\		
		\end{array}\right) $, with solution set $\xi=\left(  \begin{matrix}
			2+\oslash\\
			1+\oslash\\		
		\end{matrix}\right) $. By substitution we see that $ \xi $ also solves the system $ \left(P+(\oslash)_{2\times2}\right)x\subseteq \B $, and we 
		conclude that every real vector $ x $ satisfies $ Px\in \B $ if and only if it satisfies  $ \mathcal{R}x\subseteq \B $, with $ \mathcal{R}=\left(  \begin{array}{ll}
			1+\oslash&1+\oslash\\
			1+\oslash&-1+\oslash\\		
		\end{array}\right) $. So the systems $ P|\B $ and $ \mathcal{R}|\B $ are equivalent. The matrix $ \mathcal{R} $ is the 
		robustness matrix for the system $ P|\B $. Indeed, if a  perturbation matrix $ P+\mathcal{Q} $ of $ P $ contains a bigger neutrix than $ \oslash $, this neutrix is at least as big as $ \pounds $, and $ P+\mathcal{Q} $ is no longer limited; note also that the system $ (P+\mathcal{Q})|\B $ is inconsistent. 
		
	\end{example}
	
	We consider now reduced non-singular systems $ P|\B $. Theorem \ref{eqrepmat2} indicates that the robustness matrix $\mathcal{R}= P+E $ may be explicitly determined provided that $ \det(P) $ is not an absorber of $\overline{B}$, and the perturbations are sufficiently small such that $\mathcal{R}|\B$ remains the essential part, when adding the constraint equations generated by the neutrices of $ E $.
	\begin{theorem}\label{eqrepmat2}
		Let $ P=(a_{ij})_{n\times n} $ be a real non-singular reduced matrix and $ d=\det(P)$. Consider the system $ P|\B $ with $ \B =(\beta_1,\dots, \beta_n)^T $, where $\beta_i=b_i+B_i$ is external for all $1\leq i\leq n$ and $d$ is not an absorber of $\overline{B}$. Let $d_j=\det(M_j)$, where $M_j$ is the matrix obtained from $P$ by substituting  the $j^{th}$ column of $P$ by the column vector $b$. Let the  matrix $ E\equiv(E_{ij})_{n\times n} $ be defined by
		
		\begin{equation}\label{defEij}
			E_{ij}=\left\lbrace \begin{matrix}
				B_id/d_j&& d_j\neq 0, B_id/d_j\subset  \oslash\\
				\oslash&& \mbox{else}
			\end{matrix}\right. 
		\end{equation}	
		and $\mathcal{R}= P+E $. Assume  
		\begin{equation}\label{overB}
			B_{i}:\overline{E}_{i}\supseteq \overline{B}
		\end{equation}
		for  $ 1\leq  i\leq n $ and $d+\overline{E}$ is zeroless. Then  the non-singular matrix $\mathcal{R} $ is the  robustness matrix of $ P|\B$.
	\end{theorem}
	\begin{proof} Since $ \overline{E}\subseteq\oslash $, the matrix $ \mathcal{R}$ is reduced. So $ \det(\mathcal{R}) \subseteq d+\overline{E}$ is zeroless.
		Hence $ \mathcal{R}|\B $ is non-singular. 
		
		We show now that $ \mathcal{R}|\B $ and $ P|\B $ are equivalent. The solution of the system $P|\B$ is the external vector $ x+X $, where 
		\begin{equation}\label{xdd}
			x=P^{-1}b=(d_1/d,\dots, d_n/d)^T
		\end{equation}		
		and for $ 1\leq j\leq n $ the components $ X_{j} $ of $ X= \sum\limits_{i=1}^n B_iP^{-1}e_i$   are given by 
		\begin{equation*}\label{forX}
			X_{j}=\frac{1}{d}\sum_{k=1}^n \big(B_{k}C_{kj}\big);
		\end{equation*}
		here $C_{ij}$ is the $i, j$-cofactor of $P$. Because the cofactors of a reduced matrix are limited, and $1/d$ is not an exploder of  $\overline{B}$, it holds that 
			
		\begin{equation}\label{eqXB}
		\overline{X}\subseteq \overline{B}.
		\end{equation}
		
	 We let now $ x $ be a real vector satisfying the system $ \mathcal{R}|\B $. We define for $ 1\leq i\leq n $ a neutricial vector $ G_{i} =(G_{i1},\dots,G_{in} )^{T}$ by 
		\begin{equation*}
			G_{ij}=B_{i}:E_{ij}.
		\end{equation*}
		It follows from \eqref{defEij} that every component $ x_{j} $ of $ x $ satisfies 
		\begin{equation}\label{rescomp}
			x_{j}\in G_{ij}
		\end{equation}
		for $ 1\leq i\leq n $. Indeed, if $ d_{j}\neq 0 $ and $ B_id/d_j\subset  \oslash $ we always have $ E_{ij} = B_{i}/x_{j}$, i.e. $ x_{j}\in B_{i}:E_{ij} =G_{ij}$. If $ d_{j}\neq 0 $ and $ B_id/d_j\supseteq  \oslash $, we have $ x_{j}\oslash\subseteq B_{i}$, while $ \oslash= E_{ij}$. Again $ x_{j}\in B_{i}:E_{ij} =G_{ij}$. Finally, if  $ d_{j}= 0 $,  $ x_{j}\in G_{ij} $ is automatically satisfied, for $ x_{j}=d_{j}/d=0 $.

		Also, using \eqref{overB}, 
		\begin{equation}\label{Gmax}
			G_{ij}\supseteq B_{i}:\overline{E}_i\supseteq \overline{B} 
		\end{equation} for all $ i,j $ with $ 1\leq i,j\leq n $. 
		
		By Theorem \ref{thdec} and \eqref{rescomp} the system $ \mathcal{R}|\B $ may be decomposed into 
		
		\begin{equation}\label{sysvert}
			\left( \begin{matrix}
				P\\
				I\\
				\vdots\\
				I		
			\end{matrix} \left|  \begin{matrix}
				\mathcal B\\
				G_{1}\\
				\vdots\\
				G_{n}
			\end{matrix}\right. \right).
		\end{equation}
		In particular $ x $ satisfies $ P|\B $. Then $ P|\B $ is equivalent to $ I| (x+X)$, hence the system \eqref{sysvert} is equivalent to the system 
		\begin{equation}\label{sysI}
			\left( \begin{matrix}
				I\\
				I\\
				\vdots\\
				I		
			\end{matrix} \left|  \begin{matrix}
				x+X\\
				G_{1}\\
				\vdots\\
				G_{n}
			\end{matrix}\right. \right). 
		\end{equation}
		 If we take the representative vector $ (x,x,\dots x)^{T} $ of the right-hand side of the system \eqref{sysI}, the rank of the extended matrix of the resulting  real system is equal to $ {\rm{r}}(I|x)=n $. Clearly $ \det(I)=1 $ is not the absorber of any neutrix, and also the inclusions \eqref{Gmax} and \eqref{eqXB} hold. Then we derive from Corollary \ref{essential2} that $ I|(x+X)$ is the essential part of \eqref{sysI}. Hence \eqref{sysI} is also equivalent to $ P|\B $. We conclude that $ \mathcal{R}|\B $ and $ P|\B $ are equivalent.

		Finally we show that $\mathcal{R}$ is the robustness matrix. Let $A=(a_{ij})_{n\times n}$ be such that $(P+A)|\B$ is equivalent to $P|B$, where $ P+A $ is a limited matrix. Then for all $ i,j $ with $1\leq i, j\leq n$
		\begin{equation}\label{Aos}
			A_{ij} \subseteq \oslash.
		\end{equation}	
		The real vector $x$ given by \eqref{xdd}  is also the solution of the system   $(P+A)|\B$. Then for all $ i,j $ with $1\leq i, j\leq n$ it holds that $A_{ij}x_j=A_{ij} d_j/d\subseteq B_i$ and, if $ d_j\neq 0$, also $A_{ij}\subseteq B_i d/d_j$. Then $ A_{ij}\subseteq E_{ij} $ if  $ B_id/d_j\subset  \oslash $, and if  $ \oslash \subseteq B_id/d_j   $ or $  d_j= 0$ the inclusion	$ A_{ij}\subseteq E_{ij} $ follows from \eqref{Aos}.
		
		We conclude that $(P+E)|\B$ is the maximal limited strict perturbation of $ P|\B $, hence $\mathcal{R}$ is its robustness matrix.

	\end{proof}
	The next corollary states that the neutrices occurring in the columns of the structurally robustness matrix of a uniform system are all equal.
	\begin{corollary}\label{eqrepmat1}
		Consider the uniform system $ P|\B $, where $ P=(a_{ij})_{n\times n} $ is a real non-singular reduced matrix and  $ \B =(\beta_1,\dots, \beta_n)^T $, where $\beta_i=b_i+B$, with $ B  $ an external neutrix. Assume that $ d=\det(P)$ is not an absorber of $ B $. Let $d_i=\det(M_i)$, where $M_i$ is the matrix obtained from $P$ by substituting  the $i^{th}$ column of $P$ by the column vector $b$. Let the $n\times n$ matrix $ E $ be defined by
		$$E=\begin{pmatrix}
			E_1&\cdots & E_n\\
			\vdots &\ddots& \vdots\\
			E_1&\cdots& E_n
		\end{pmatrix},$$
		where $ E_i=\min(\oslash, Bd/d_i) $ for $1\leq i\leq n$. Let $\overline{E}=\max\limits_{1\leq i\leq n} E_i$ and $\mathcal{R}= P+E $. Assume $d+\overline{E}$ is zeroless. Then the non-singular matrix $\mathcal{R} $ is the robustness matrix of $ P|\B $.
	\end{corollary} 
	\begin{proof} The result follows from Theorem \ref{eqrepmat2}, observing that 
		$E_{1j}=E_{2j}=\cdots=E_{nj}=E_j$ for all $1\leq j\leq n$, and that condition \eqref{overB} is satisfied, since $	B:\overline{E}_i =
		B: \overline{E}	\supseteq B:\oslash\supseteq B$ for $1\leq  i\leq n$. 
	\end{proof}
	We end with some examples. We start with a flexible system for which the robustness matrix $\mathcal{R} $ exhibits different neutrices for each entry. Then we give an example of a uniform system with a robustness matrix having identical neutrices in each column. The final example shows that Theorem \ref{eqrepmat2} is no longer valid if the determinant  of the coefficient matrix is an absorber of the neutrices at the right-hand side of the flexible system. 
	\begin{example}	
		Consider the system 
		\begin{equation*}\label{omega22N}
			P|\B=
			\left( \begin{array}{ll|l}
				1&1&\omega+2+\pounds/\omega\\
				1&-1&\omega+\oslash\\		
			\end{array}\right). 
		\end{equation*}
		Then \begin{equation*}
			\mathcal{R}\equiv
			\left( \begin{array}{ll}
				1+\pounds/\omega^{2}&1+\pounds/\omega\\
				1+\oslash/\omega&-1+\oslash\\		
			\end{array}\right) 
		\end{equation*}
		is the robustness matrix of \(Q|\B\). Indeed $\det(P)=-2 $, which is not an absorber of any neutrix and the verification of the remaining conditions of Theorem \ref{eqrepmat2} is straightforward.
	\end{example}

	\begin{example} Let $ \omega\in \R^{+} $ be unlimited. Consider the uniform system \begin{equation*}\label{omega22}
			P|\B=
			\left(  \begin{array}{ll|l}
				1&1&\omega+2+\oslash\\
				1&-1&\omega+\oslash\\		
			\end{array}\right).
		\end{equation*}
		It is straightforward to verify that by Corollary \ref{eqrepmat1} we obtain the robustness matrix
		\begin{equation*}
			\mathcal{R}=
			\left(  \begin{array}{ll}
				1+\oslash/\omega&1+\oslash\\
				1+\oslash/\omega&-1+\oslash\\		
			\end{array}\right).
		\end{equation*}

	\end{example}
	
	\begin{example}	Let \(\eps\) be a positive infinitesimal. Consider the uniform flexible system 
		\begin{equation*}P|\B=
			\left( 	\begin{array}{ll|l}
				1&1 &1+\oslash\\
				1&1+\eps& 1+\oslash
			\end{array} \right).
		\end{equation*}
		We have \(d\equiv\det\begin{pmatrix}
			1&1\\ 1& 1+\eps
		\end{pmatrix}=\eps\neq 0\). Let us choose \(b=\begin{pmatrix}
			1\\1
		\end{pmatrix}\) as a representative of the right-hand side. Then applying Cramer's rule we obtain  \(d_1=\eps\) and \(d_2=0\). Suppose we define \( E=(E_{ij})_{1\leq i,j\leq 2} \) as in \eqref{defEij}, then \(E_{ij}=\oslash\) for \(1\leq i, j\leq 2\). he fact that $d$ is an absorber of the neutrix $ \oslash $ at the right-hand side and the fact that \(\det(P+E)=\oslash\) imply that two conditions of Corollary \ref{eqrepmat1} are not satisfied. In addition, it is obvious that the matrix $ P+E $ is singular. We show that \begin{equation*}		
			(P+E)|\B=\left(  \begin{array}{ll|l}
				1+\oslash &1+\oslash &1+\oslash\\
				1+\oslash &1+\oslash &1+\oslash\\
			\end{array} \right)	\end{equation*} is not equivalent to $ P|\B $. 
			Indeed, let $ \xi_{P} $ is the solution of $ P|\B $ and $ \xi_{P+E} $ is the solution of $ (P+E)|\B $. A straigthforward application of the parameter method of Theorem \ref{thesolhom} shows that  \(\xi_{P} =\begin{pmatrix}
			1\\0
		\end{pmatrix}+\oslash/\eps \begin{pmatrix}
			1+\eps\\-1
		\end{pmatrix}+\oslash/\eps  \begin{pmatrix}
			-1\\1
		\end{pmatrix}\), from which we derive that \(\xi_{P} =\begin{pmatrix}
			1\\0
		\end{pmatrix}+\oslash \begin{pmatrix}
			1\\0
		\end{pmatrix}+\oslash/\eps \begin{pmatrix}
			-1\\1
		\end{pmatrix}\). Also, the singular system $ (P+E)|\B $ is equivalent with 
	\begin{equation*}		
			P'|\B'\equiv\left( \begin{array}{ll|l}
				1+\oslash &1+\oslash &1+\oslash
			\end{array} \right).	\end{equation*}
		Its associated integrated system becomes		
		 \begin{equation*}		
			\left( \begin{matrix}
				P'\\
				K 
			\end{matrix}\left|   \begin{matrix}
				\mathcal B'\\
				\F^{(c)}
			\end{matrix}\right.\right) =\left(  \begin{array}{ll|l}
				1 &1 &1+\oslash\\
				1 &0 &\pounds\\
				0 &1 &\pounds\\
			\end{array} \right).	\end{equation*}
		If we put it in increasing row-echelon form and apply the parameter method we find that 
		\(\xi_{P+E}=\begin{pmatrix}
			1\\0
		\end{pmatrix}+\oslash \begin{pmatrix}
			1\\0
		\end{pmatrix}+\pounds \begin{pmatrix}
			-1\\1
		\end{pmatrix}\).
		
		We see that $ \xi_{P+E} \subset  \xi_{P} $, hence $ P|\B $ and $(P+E)|\B  $ are not equivalent.
	\end{example}
	
\end{document}